\documentclass[a4paper,12pt]{amsart}

\usepackage{float}
\usepackage{euscript,eufrak,verbatim}
\usepackage{graphicx}
\usepackage[usenames]{color}
\usepackage[colorlinks,linkcolor=red,anchorcolor=blue,citecolor=blue]{hyperref}
\usepackage{amsmath}
\usepackage{amsthm}
\usepackage[all]{xy}
\usepackage{graphicx} 
\usepackage{amssymb} 

\usepackage{mathrsfs}
\usepackage{amscd}

\theoremstyle{plain}
\newtheorem{main theorem}{Main Theorem}
\newtheorem{theorem}{Theorem}[section]
\newtheorem{lemma}[theorem]{Lemma}

\newtheorem{corollary}[theorem]{Corollary}
\newtheorem{proposition}[theorem]{Proposition}
\newtheorem{claim}[theorem]{Claim}

\theoremstyle{definition}

\newtheorem{remark}[theorem]{Remark}
\newtheorem{example}[theorem]{Example}

\makeatletter 
\@addtoreset{equation}{section}
\numberwithin{equation}{section}

\newcommand{\norm}[1]{\left\lVert#1\right\rVert}
\newcommand{\diam}{\mathrm{Diam}}

\newcommand{\mdim}{\mathrm{mdim}}

\newcommand{\mdimh}{\mathrm{mdim}_{\mathrm{H}}}
\newcommand{\mdimm}{\mathrm{mdim}_{\mathrm{M}}}
\newcommand{\htop}{h_{\mathrm{top}}}
\newcommand{\dimh}{\mathrm{dim}_{\mathrm{H}}}
\newcommand{\dimm}{\mathrm{dim}_{\mathrm{M}}}
\newcommand{\umdimm}{\overline{\mathrm{mdim}}_{\mathrm{M}}}
\newcommand{\lmdimm}{\underline{\mathrm{mdim}}_{\mathrm{M}}}
\newcommand{\umdimh}{\overline{\mathrm{mdim}}_{\mathrm{H}}}
\newcommand{\lmdimh}{\underline{\mathrm{mdim}}_{\mathrm{H}}}
\newcommand{\dist}{\mathrm{dist}}

\title[Infinite dimensional fractals]{Mean Hausdorff dimension of some infinite dimensional fractals}

\author{Masaki Tsukamoto}

\address
{Department of Mathematics, Kyoto University, Kitashirakawa Oiwake-cho, Sakyo-ku, Kyoto 606-8502, Japan}

\email{tsukamoto@math.kyoto-u.ac.jp}

\begin{document}



\thanks{M.T. was supported by JSPS KAKENHI JP21K03227.}

\maketitle

\begin{abstract}
Mean Hausdorff dimension is a dynamical version of Hausdorff dimension.
It provides a way to dynamicalize geometric measure theory.
We pick up the following three classical results of fractal geometry.
 \begin{enumerate}
    \item The calculation of Hausdorff dimension of homogeneous sets in the circle.
    \item The coincidence of Hausdorff and Minkowski dimensions for self-similar sets.
    \item The calculation of Hausdorff dimension of Bedford--McMullen carpets. 
 \end{enumerate}
We develop their analogues for mean Hausdorff dimension:
 \begin{enumerate}
    \item[(1’)] The calculation of mean Hausdorff dimension of homogeneous systems in the infinite dimensional torus.
    \item[(2’)] The coincidence of mean Hausdorff dimension and metric mean dimension for self-similar systems.
    \item[(3’)] The calculation of mean Hausdorff dimension of infinite dimensional carpets.
 \end{enumerate}
\end{abstract}

\section{Background} \label{section: background}

In the classical dimension theory, there are three famous notions of dimension:
topological dimension, Hausdorff dimension and Minkowski dimension.
Their basic relation is the following.
\[  \text{Topological dimension} \leq \text{Hausdorff dimension} \leq 
     \text{Minkowski dimension}. \]

At the end of the 20th century, Gromov \cite{Gromov} found a way to \textit{dynamicalize}
topological dimension theory.
He introduced a dynamical version of topological dimension called \textit{mean dimension}.
Given a dynamical system, its mean dimension measures how many parameters per iterate
we need for describing its orbits.
Mean dimension is also called mean topological dimension.

Lindenstrauss--Weiss \cite{Lindenstrauss--Weiss} introduced a dynamical version of 
Minkowski dimension called \textit{metric mean dimension}
in order to better understand relations between mean dimension and topological entropy.
The definition of metric mean dimension is a fusion of the definitions of topological entropy and Minkowski dimension.
We will review it in \S \ref{section: basic definitions}.

In order to connect mean dimension to ergodic theory, 
Lindenstrauss--Tsukamoto \cite{Lindenstrauss--Tsukamoto_double} introduced a dynamical version of 
Hausdorff dimension called \textit{mean Hausdorff dimension}.
See \S \ref{section: basic definitions} for the definition. 
Mean Hausdorff dimension is better suited for measure theoretic studies than metric mean dimension, 
as (ordinary) Hausdorff dimension is more closely related to measure theory than Minkowski dimension.

The following is the basic relation between the above three dynamical dimensions.
\[  \text{mean dimension} \leq \text{mean Hausdorff dimension} \leq 
     \text{metric mean dimension}. \]
See Proposition \ref{prop: mean Hausdorff dimension and metric mean dimension} and 
Remark \ref{remark: mean dimension and mean Hausdorff dimension}
in \S \ref{section: basic definitions} for more precise statements.
     
We mainly study mean Hausdorff dimension and metric mean dimension in this paper.
Mean Hausdorff dimension is more difficult to evaluate than metric mean dimension.
So we often pay more attention to mean Hausdorff dimension.

We pick up the following three classical results of fractal geometry.

\begin{itemize}
  \item \textbf{Homogeneous sets} (\cite{Furstenberg}).
   Let $\mathbb{R}/\mathbb{Z}$ be the circle with a metric 
   \[  d(x, y) := \min_{n\in \mathbb{Z}} |x-y-n|. \]
   Let $a>1$ be a natural number greater than one. We consider the “$\times a$ map” on the circle $\mathbb{R}/\mathbb{Z}$:
   \[  T_a:\mathbb{R}/\mathbb{Z} \to \mathbb{R}/\mathbb{Z}, \quad 
        x\mapsto ax. \]
    Suppose $X\subset \mathbb{R}/\mathbb{Z}$ is a closed subset satisfying $T_a(X)\subset X$.
   Furstenberg \cite[Proposition III.1]{Furstenberg} proved that 
   its Hausdorff dimension $\dim_{\mathrm{H}}(X,d)$ coincides with the Minkowski dimension $\dim_{\mathrm{M}}(X,d)$
   and that they are  given by 
   \[  \dim_\mathrm{H}(X, d) = \dim_{\mathrm{M}}(X,d)  = \frac{\htop(X,T_a)}{\log a}. \]
   Here $\htop(X,T_a)$ is the topological entropy of the dynamical system $(X, T_a)$.

  \item \textbf{Self-similar sets} (\cite{Falconer}). 
  Let $f_1, f_2, \dots, f_k$ be contracting similarity transformations of the Euclidean space $\mathbb{R}^n$.
  Let $X$ be the attractor of the family $\{f_1, f_2, \dots, f_n\}$.
  Namely, $X$ is the unique nonempty and compact subset of $\mathbb{R}^n$
  satisfying 
  \[  X = \bigcup_{i=1}^k f_i(X). \]
   Let $d$ be the Euclidean metric on $\mathbb{R}^n$. 
   Falconer \cite[Example 2]{Falconer} proved that the Hausdorff dimension of 
  $(X, d)$ is equal to its Minkowski dimension:
   \[  \dim_{\mathrm{H}}(X,d) = \dim_{\mathrm{M}}(X,d). \]

  \item \textbf{Bedford--McMullen carpets} (\cite{Bedford, McMullen}).
  Let $a$ and $b$ be natural numbers with $a\geq b\geq 2$.
  Set $A = \{0, 1, 2, \dots, a-1\}$ and $B = \{0,1,2,\dots, b-1\}$.
  Let $R\subset A\times B$ be a non-empty subset. 
 We define a closed subset $X$ of the unit square $[0,1]^2$ by 
  \[  X = \left\{\left(\sum_{n=1}^\infty \frac{x_n}{a^n}, \sum_{n=1}^\infty \frac{y_n}{b^n}\right) \in [0,1]^2 \, \middle|\, \forall n\geq 1:
             (x_n, y_n)\in R\right\}. \]
  This space was first introduced and studied by Bedford \cite[Chapter 4]{Bedford} and McMullen \cite{McMullen}, so it is called 
  a \textit{Bedford--McMullen carpet}.
  It is a famous example of fractal sets whose Hausdorff and Minkowski dimensions do not coincide.
  Let $d$ be the Euclidean metric on the plane.
  For each $j\in B$ we denote by $t_j$ the number of $i\in A$ with $(i, j)\in R$.
  The Hausdorff dimension of $(X, d)$ is given by 
  \[  \dim_{\mathrm{H}}(X,d) = \log_b\left(\sum_{j=0}^{b-1} t_j^{\log_a b}\right). \]
  Let $r$ be the cardinality of $R$, and let $s$ be the number of $j \in B$ for 
  which there exists $i\in A$ with $(i, j)\in R$.
  The Minkowski dimension is given by 
  \[  \dim_{\mathrm{M}}(X, d) = \log_b s + \log_a \left(\frac{r}{s}\right). \]
  Except for some special cases\footnote{Namely, either $n=m$ or all nonzero $t_j$ are equal to each other.}, 
  the Hausdorff dimension $\dim_{\mathrm{H}}(X, d)$ 
  is strictly smaller than the Minkowski dimension $\dim_{\mathrm{M}}(X, d)$.
\end{itemize}

The purpose of this paper is to develop analogues of these results for mean Hausdorff dimension and metric mean dimension.
Our main results are Theorem 
\ref{theorem: mean Hausdorff dimension of homogeneous systems} 
(an analogue of Furstenberg's theorem),
Theorem \ref{main theorem for self-similar systems} 
(an analogue of Falconer's theorem) and 
Theorem 
\ref{theorem: mean Hausdorff dimension and metric mean dimension of carpet systems}
(the mean Hausdorff dimension of ``infinite dimensional carpets'')
below.

This paper is a starting point for the study of “infinite dimensional fractals”.
Our primary purpose is just to show how to formulate meaningful mathematical theorems
about infinite dimensional fractals.
Hopefully such study will become a fruitful research area in a future.

The plan of this paper is as follows:
In \S \ref{section: basic definitions} we explain the definitions of mean Hausdorff dimension and metric mean dimension.
In \S \ref{section: mean Hausdorff dimension of homogeneous systems} we study an analogue of Furstenberg’s theorem for mean 
Hausdorff dimension.
In \S \ref{section: self-similarity and mean Hausdorff dimension} we introduce a dynamical version of self-similar sets
and study an analogue of Falconer’s theorem.
In \S \ref{section: infinite dimensional carpets} 
we study the mean Hausdorff dimension of infinite dimensional carpets.
In the Appendix we explain one more example of the calculations of mean Hausdorff dimension.

\section{Basic definitions}  \label{section: basic definitions}

The purpose of this section is to review the definitions of metric mean dimension (\cite{Lindenstrauss--Weiss})
and mean Hausdorff dimension (\cite{Lindenstrauss--Tsukamoto_double}).

First we prepare some basic quantities of compact metric spaces.
Let $(X, d)$ be a compact metric space.
For a positive number $\varepsilon$ we define the \textbf{$\varepsilon$-covering number} 
$\#(X,d,\varepsilon)$ as the minimum cardinality $n$ of open covers
$\{U_1, \dots, U_n\}$ of $X$ satisfying $\diam\, U_i < \varepsilon$ for all $1\leq i\leq n$. 
We define the \textbf{$\varepsilon$-scale Minkowski dimension of $X$} by 
\[  \dimm(X,d,\varepsilon) = \frac{\log \#(X, d, \varepsilon)}{\log(1/\varepsilon)}. \]
The \textbf{upper and lower Minkowski dimensions} of $(X,d)$ are defined by 
\[ \overline{\dim}_{\mathrm{M}}(X,d) = \limsup_{\varepsilon\to 0} \dimm(X, d, \varepsilon), \quad 
    \underline{\dim}_{\mathrm{M}}(X,d) = \liminf_{\varepsilon\to 0} \dimm(X, d, \varepsilon).  \]
When these two values coincide, it is called the \textbf{Minkowski dimension} of $(X,d)$ and denoted by $\dim_{\mathrm{M}}(X,d)$.

For $s\geq 0$ and $\varepsilon>0$
we define $\mathcal{H}_\varepsilon^s(X,d)$ by 
\[ \mathcal{H}^s_\varepsilon(X, d) = \inf\left\{\sum_{i=1}^\infty \left(\diam\, E_i\right)^s\, \middle|\, 
                                               X = \bigcup_{i=1}^\infty E_i \text{ with } \diam \, E_i < \varepsilon\> (\forall i\geq 1)\right\}. \]
The meaning of the term $\left(\diam\, E_i\right)^s$ becomes ambiguous when $s=0$ (with $\diam \, E_i = 0$) or $E_i = \emptyset$.
We use the convention that $0^0 = 1$ and $\left(\diam\, \emptyset\right)^s  = 0$ for all $s\geq 0$.
Note that this convention implies $\mathcal{H}^0_\varepsilon(X, d)\geq 1$ (assuming $X\neq \emptyset$).
For $\varepsilon>0$ we define the 
\textbf{$\varepsilon$-scale Hausdorff dimension} $\dimh(X,d,\varepsilon)$ as the supremum of $s\geq 0$ 
satisfying $\mathcal{H}^s_\varepsilon(X, d) \geq 1$.
The \textbf{Hausdorff dimension} of $(X,d)$ is defined by 
\[  \dimh(X,d) = \lim_{\varepsilon\to 0} \dimh(X,d,\varepsilon). \]
For $0<\varepsilon <1$ we have 
\begin{equation} \label{eq: Hausdorff and Minkowski}
    \dimh(X,d,\varepsilon) \leq \dimm(X, d, \varepsilon).
\end{equation}
Hence 
\[  \dimh(X,d) \leq \underline{\dim}_{\mathrm{M}}(X,d) \leq \overline{\dim}_{\mathrm{M}}(X,d). \]

Next we consider dynamical versions of Hausdorff and Minkowski dimensions.
A pair $(X, T)$ is called a \textbf{dynamical system} if $X$ is a compact metrizable space and 
$T:X\to X$ is a continuous map\footnote{In many literature of mean dimension theory 
(e.g. \cite{Lindenstrauss--Weiss, Lindenstrauss, Lindenstrauss--Tsukamoto_double}), 
one usually considers only invertible dynamical systems 
(namely, the case that $T$ is a homeomorphism).
But in this paper we do not assume that $T$ is a homeomorphism.
This is mainly because the “$\times a$ map” in Furstenberg’s theorem 
is non-invertible and 
we would like to study its variation in mean dimension theory.}.

Let $(X, T)$ be a dynamical system with a metric $d$ on $X$.
For each natural number $N$ we define a metric $d_N$ on $X$ by 
\[  d_N(x, y) = \max_{0\leq n <N} d\left(T^n x, T^n y\right). \]
We sometimes use the notation $d^T_N$ for $d_N$ in order to clarify the map $T$.
The \textbf{topological entropy} of $(X, T)$ is defined by 
\[  \htop(X,T) = \lim_{\varepsilon\to 0} \left(\lim_{N\to \infty} \frac{\log \#\left(X,d_N,\varepsilon\right)}{N}\right). \]
The topological entropy is a topological invariant of dynamical systems.
Namely the value is independent of the choice of a metric $d$.

We define the \textbf{upper and lower metric mean dimensions} of $(X, T, d)$ by
\begin{align*}
     \umdimm(X,T,d) & = \limsup_{\varepsilon\to 0} 
     \left(\lim_{N\to \infty} \frac{\dimm\left(X,d_N,\varepsilon\right)}{N}\right), \\
     \lmdimm(X,T,d) & = \liminf_{\varepsilon\to 0} 
     \left(\lim_{N\to \infty} \frac{\dimm\left(X,d_N,\varepsilon\right)}{N}\right).
\end{align*}
These values depend on the choice of a metric $d$.
When the upper and lower metric mean dimensions coincide, the common value is called the \textbf{metric mean dimension} of $(X, T, d)$ and denoted by 
$\mdimm(X, T, d)$.

We define the \textbf{upper and lower mean Hausdorff dimensions} of $(X, T, d)$ by 
\begin{align*}
   \umdimh(X, T, d) = \lim_{\varepsilon\to 0} \left(\limsup_{N\to \infty} \frac{\dimh\left(X, d_N,\varepsilon\right)}{N}\right), \\
   \lmdimh(X, T, d) = \lim_{\varepsilon\to 0} \left(\liminf_{N\to \infty} \frac{\dimh\left(X, d_N,\varepsilon\right)}{N}\right).
\end{align*}
These also depend on the choice of $d$.
When they coincide, the common value is called the \textbf{mean Hausdorff dimension} of $(X,T,d)$ and denoted by 
$\mdimh(X, T, d)$.

\begin{proposition} \label{prop: mean Hausdorff dimension and metric mean dimension}
\[   \lmdimh(X, T, d) \leq \umdimh(X, T, d) \leq  \lmdimm(X,T,d)  \leq \umdimm(X,T,d). \]
\end{proposition}

\begin{proof}
Let $0<\varepsilon<1$.
From (\ref{eq: Hausdorff and Minkowski}) 
\[   \frac{\dimh(X,d_N,\varepsilon)}{N}  \leq \frac{\dimm(X,d_N,\varepsilon)}{N}.  \]
Hence we have $\umdimh(X, T, d) \leq  \lmdimm(X,T,d)$.
The rests are trivial.
\end{proof}

\begin{remark} \label{remark: mean dimension and mean Hausdorff dimension}
We denote by $\mdim(X, T)$ the mean topological dimension of a dynamical system $(X, T)$.
Then for any metric $d$ on $X$ we have \cite[Proposition 3.2]{Lindenstrauss--Tsukamoto_double}
\[   \mdim(X, T) \leq \lmdimh(X, T, d). \]
We do not explain the definition of mean topological dimension here because we will not use it in the sequel.
\end{remark}

Throughout the paper we denote the set of natural numbers by $\mathbb{N} = \{1,2,3,\dots\}$.

\begin{example} \label{example: Hilbert cube}
Let $[0,1]^\mathbb{N} = [0,1]\times [0,1]\times [0,1]\times \cdots$ be the infinite dimensional cube.
We define the shift map $\sigma$ on it by 
\[  \sigma\left((x_n)_{n\in \mathbb{N}}\right) = (x_{n+1})_{n\in \mathbb{N}}. \]
We define a metric $d$ on $[0,1]^\mathbb{N}$ by 
\[  d\left((x_n)_{n\in \mathbb{N}}, (y_n)_{n\in \mathbb{N}}\right) = \sum_{n=1}^\infty 2^{-n}|x_n-y_n|. \]
Then 
\[  \mdim\left([0,1]^\mathbb{N}, \sigma\right) = \mdimh\left([0,1]^\mathbb{N},\sigma, d\right) = \mdimm\left([0,1]^\mathbb{N},\sigma,d\right) = 1. \]
\end{example}

\begin{example} \label{example: 0,1,1/2,1/3,dots}
Let $K = \{0\}\cup\{\frac{1}{n}\mid n\geq 1\} = \{0,1,\frac{1}{2},\frac{1}{3},\dots\}$.
It is well-known that the Hausdorff dimension and Minkowski dimension of $K$ with respect to the Euclidean metric are 
zero and $\frac{1}{2}$ respectively.
We consider 
\[  K^\mathbb{N} = K\times K\times K \times \cdots. \]
We define the shift map $\sigma: K^{\mathbb{N}} \to K^{\mathbb{N}}$ and a metric $d$ on $K^\mathbb{N}$ 
as in Example \ref{example: Hilbert cube}.
Then 
\[  \mdimh\left(K^\mathbb{N},\sigma,d\right) = 0, \quad
     \mdimm\left(K^\mathbb{N},\sigma,d\right) = \frac{1}{2}. \]
The proof of $\mdimm\left(K^\mathbb{N},\sigma,d\right) = \frac{1}{2}$ is straightforward.
But it is not so easy to prove $\mdimh\left(K^\mathbb{N},\sigma,d\right) = 0$.
We will explain it in the Appendix.
\end{example}

\section{Mean Hausdorff dimension of homogeneous systems}  \label{section: mean Hausdorff dimension of homogeneous systems}

We develop an analogue of Furstenberg’s theorem \cite[Proposition III.1]{Furstenberg} in this section.

\subsection{Topological entropy of $\mathbb{N}^2$-actions} \label{subsection: topological entropy of N^2-actions}

We need to introduce topological entropy of $\mathbb{N}^2$-actions in order to explain the main result of this section.
A triple $(X,S, T)$ is called a \textbf{$\mathbb{N}^2$-action} if $X$ is a compact metrizable space, and if
$S:X\to X$ and $T:X\to X$ are continuous maps with $S\circ T = T\circ S$.

Let $(X, S, T)$ be a $\mathbb{N}^2$-action with a metric $d$ on $X$.
For a subset $\Omega\subset \mathbb{N}^2$ we define a metric $d^{S, T}_\Omega$ on $X$ by 
\[  d^{S,T}_\Omega(x,y) = \sup_{(m,n)\in \Omega} d\left(S^m T^n x, S^m T^n y\right). \]
It is convenient to use this notation also for the case that $\Omega$ is a subset of $\mathbb{R}^2$:
For a subset $\Omega\subset \mathbb{R}^2$ we set
\[  d^{S,T}_\Omega(x,y) 
     = \sup_{(m,n)\in \Omega\cap \mathbb{N}^2} d\left(S^m T^n x, S^m T^n y\right). \]

We define the \textbf{topological entropy of $(X, S, T)$} by 
\[  \htop(X, S, T) = \lim_{\varepsilon\to 0} \left(\lim_{N\to \infty}\frac{\log\#\left(X,d^{S,T}_{[0,N)^2}, \varepsilon\right)}{N^2}\right). \]
It is easy to check that the limits exist. 
We also have 
\[ \htop(X, S, T) = \lim_{\varepsilon\to 0}
                          \left(\lim_{\substack{M\to \infty \\ N\to \infty}}
                          \frac{\log \#\left(X, d^{S,T}_{[0,M)\times [0,N)},\varepsilon\right)}{MN}
                          \right). \]
Here 
\begin{align*}
   d^{S,T}_{[0,N)^2}(x, y) &= \max_{\substack{0\leq m <N \\ 0\leq n< N}} 
    d\left(S^m T^n x, S^m T^n y\right),  \\ 
    d^{S,T}_{[0,M)\times [0,N)}(x,y) &= \max_{\substack{0\leq m <M \\ 0\leq n <N}}
         d\left(S^m T^n x, S^m T^n y\right).  
\end{align*}

\subsection{Main result for homogeneous systems}  \label{subsection: main result of homogeneous systems}
Let $\mathbb{R}/\mathbb{Z}$ be the circle with a metric $\rho$ defined by 
\[  \rho(x, y) = \min_{n\in \mathbb{Z}} |x-y-n|. \]
We consider the “infinite dimensional torus”: 
\[   \left(\mathbb{R}/\mathbb{Z}\right)^{\mathbb{N}}
    = \mathbb{R}/\mathbb{Z}\times \mathbb{R}/\mathbb{Z}\times \mathbb{R}/\mathbb{Z}\times \cdots. \]
We define a metric $d$ on it by 
\[   d\left((x_n)_{n\in \mathbb{N}}, (y_n)_{n\in \mathbb{N}}\right) = \sum_{n=1}^\infty 2^{-n} \rho(x_n, y_n). \]
We define the shift map $\sigma: \left(\mathbb{R}/\mathbb{Z}\right)^\mathbb{N}\to  \left(\mathbb{R}/\mathbb{Z}\right)^\mathbb{N}$ by
\[  \sigma\left((x_n)_{n\in \mathbb{N}}\right) = (x_{n+1})_{n\in \mathbb{N}}. \]

Let $a>1$ be a natural number greater than one.
We define the “$\times a$ map” $T_a: \left(\mathbb{R}/\mathbb{Z}\right)^\mathbb{N} \to \left(\mathbb{R}/\mathbb{Z}\right)^\mathbb{N}$ by 
\[  T_a\left((x_n)_{n\in \mathbb{N}}\right) = (ax_n)_{n\in \mathbb{N}}. \]
Notice that $T_a$ and $\sigma$ commute.

The following is the main result of this section.

\begin{theorem} \label{theorem: mean Hausdorff dimension of homogeneous systems}
Let $X\subset  \left(\mathbb{R}/\mathbb{Z}\right)^\mathbb{N}$ be a closed subset such that 
$\sigma(X) \subset X$ and $T_a(X)\subset X$. Then 
\[  \mdimh(X,\sigma,d) = \mdimm(X,\sigma,d) = \frac{\htop(X, \sigma,T_a)}{\log a}. \]
Here $\htop(X, \sigma,T_a)$ is the topological entropy of the $\mathbb{N}^2$-action $(X,\sigma, T_a)$.
\end{theorem}

A “symbolic dynamics version” of this theorem was presented in \cite{Shinoda--Tsukamoto}.

The proof of Theorem \ref{theorem: mean Hausdorff dimension of homogeneous systems} consists of two parts:
the proofs of the upper bound 
\begin{equation}  \label{eq: upper bound on metric mean dimension of homogeneous system}
     \umdimm(X, T, d) \leq \frac{\htop(X, \sigma,T_a)}{\log a}   
\end{equation}     
and the lower bound 
\begin{equation}  \label{eq: lower bound on mean Hausdorff dimension of homogeneous system}
     \lmdimh(X, T, d) \geq  \frac{\htop(X, \sigma,T_a)}{\log a}. 
\end{equation}     

The upper bound (\ref{eq: upper bound on metric mean dimension of homogeneous system}) directly follows from the definitions.
The proof of the lower bound (\ref{eq: lower bound on mean Hausdorff dimension of homogeneous system}) is more involved.
The next subsection is a preparation for it.

\subsection{Lipschitz map and mean Hausdorff dimension}  
\label{subsection: Lipschitz map and mean Hausdorff dimension}

Dai--Zhou--Geng \cite[Theorem 2]{Dai--Zhou--Geng} proved that if $T:X\to X$ is a Lipschitz map 
of a compact metric space $(X, d)$ with a Lipschitz constant $L>1$ then
\begin{equation} \label{eq: entropy is bounded by hausdorff dimension}
   \dimh(X,d) \geq \frac{\htop(X, T)}{\log L}. 
\end{equation}   
See also \cite[Corollary 2.2]{Misiurewicz}.
The purpose of this subsection is to prove a variation of this result for 
mean Hausdorff dimension.

First we prove a “finite accuracy version” of (\ref{eq: entropy is bounded by hausdorff dimension}).

\begin{lemma}  \label{lemma: entropy is bounded by hausdorff dimension finitely version}
Let $(X, T)$ be a dynamical system with a metric $d$ on $X$.
Suppose there is $L>1$ satisfying 
\[  d(T x, T y) \leq L \cdot d(x, y), \quad (x, y\in X). \]
Let $t, \delta, \varepsilon$ be positive numbers satisfying 
\[  0< t<1, \quad 0<\delta < 1, \quad \delta^{1-t} < \varepsilon. \]
Then 
\begin{equation} \label{eq: entropy is bounded by hausdorff dimension finitely version}
    \lim_{N\to \infty} \frac{\log \#\left(X, d^T_N, \varepsilon\right)}{N} \leq 
      \frac{\log L}{t} \cdot \dimh(X, d, \delta). 
\end{equation}      
\end{lemma}

Notice that if we let $\delta\to 0$ in the inequality 
(\ref{eq: entropy is bounded by hausdorff dimension finitely version}) then we get
\[   \lim_{N\to \infty} \frac{\log \#\left(X, d^T_N, \varepsilon\right)}{N}
     \leq \frac{\log L}{t} \cdot \dimh(X, d). \]
Letting $t\to 1$ and $\varepsilon\to 0$, we get 
(\ref{eq: entropy is bounded by hausdorff dimension}).
So (\ref{eq: entropy is bounded by hausdorff dimension}) follows from 
(\ref{eq: entropy is bounded by hausdorff dimension finitely version}).

\begin{proof}[Proof of Lemma \ref{lemma: entropy is bounded by hausdorff dimension finitely version}]
The following proof is motivated by the arguments of 
\cite[Proposition III.1]{Furstenberg} and 
\cite[Proposition 1]{Bowen}.

Let $s$ be a positive number satisfying $\dimh(X, d, \delta)<s$.
There exists an open cover $X = U_1\cup \dots \cup U_M$ such that 
$\diam(U_m,d)<\delta$ for all $1\leq m\leq M$ and
\[  \sum_{m=1}^M \left(\diam(U_m, d)\right)^s < 1. \]
Choose positive numbers $\delta_m$ $(1\leq m\leq M)$ such that $\diam(U_m, d) < \delta_m < \delta$ and 
\[  \sum_{m=1}^M \delta_m^s < 1. \]
Set 
\[  N_m = \lceil \log_L \delta_m^{-t} \rceil \quad (\geq 1). \]
Here $\lceil x\rceil = \min\{n\in \mathbb{Z}\mid n\geq x\}$ for real numbers $x$.
Then $L^{N_m} \geq \delta_m^{-t}$ and hence
\begin{equation}  \label{eq: N_m are large}
   \sum_{m=1}^M L^{-sN_m/t} \leq  \sum_{m=1}^M \delta_m^s < 1.
\end{equation}
From the Lipschitz condition of $T$
\begin{align*}
   \diam\left(U_m,d^T_{N_m}\right) &\leq L^{N_m-1} \cdot \diam(U_m,d) \\
    &\leq L^{\log_L \delta_m^{-t}} \cdot \diam(U_m,d)  \quad
    (\text{by } N_m= \lceil \log_L \delta_m^{-t} \rceil \leq \log_L \delta_m^{-t}+1) \\
    &= \delta_m^{-t} \cdot \diam(U_m,d) \\
    &< \delta_m^{1-t} < \delta^{1-t} < \varepsilon.
\end{align*}

Let $N$ be a natural number.
Let $I_N$ be the set of sequences $(m_1, m_2, \dots, m_k)$ of 
natural numbers $m_i$ satisfying 
\begin{itemize}
   \item $1\leq m_i \leq M$ for all $1\leq i \leq k$,
   \item $N_{m_1}+N_{m_2}+\dots + N_{m_{k-1}} < N \leq N_{m_1}+N_{m_2}+\dots + N_{m_{k-1}} + N_{m_k}$.
\end{itemize}
Here $k$ is not a fixed number. It also varies.

We define an open covering $\mathcal{U}$ of $X$ as the family of
\[   U_{m_1}\cap T^{-N_{m_1}}U_{m_2}\cap T^{-N_{m_1}-N_{m_2}}U_{m_3}\cap \dots \cap T^{-N_{m_1}-N_{m_2}-\dots-N_{m_{k-1}}}U_{m_k},  \]
where $(m_1, m_2, \dots, m_k)\in I_N$.
Every member $U\in \mathcal{U}$ satisfies $\diam\left(U, d^T_N\right) < \varepsilon$.
Hence $\#\left(X, d_N^T, \varepsilon\right) \leq \left|I_N\right|$.
Here $|I_N|$ denotes the cardinality of $I_N$.

Set $\bar{N} := \max\left(N_1,\dots,N_M\right)$.
We have 
\begin{align*}
  |I_N| \cdot L^{-s(N+\bar{N})/t}  & \leq \sum_{(m_1,\dots,m_k)\in I_N} 
  L^{-s(N_{m_1}+\dots+N_{m_k})/t} \\
   &\leq \sum_{k=1}^\infty \left(\sum_{m=1}^\infty L^{-s N_m/t}\right)^k < \infty \quad 
   (\text{by (\ref{eq: N_m are large})}). 
\end{align*}  
Hence 
\[  \log |I_N| \leq \frac{s(N+\bar{N})}{t} \cdot \log L + \mathrm{const},  \]
where $\mathrm{const}$ is a positive constant independent of $N$.
Therefore 
\[  \log \#\left(X, d^T_N, \varepsilon \right) \leq 
     \frac{s(N+\bar{N})}{t} \cdot \log L + \mathrm{const}.  \]
Dividing this by $N$ and letting $N\to \infty$, we get 
\[  \lim_{N\to \infty} \frac{\log \#\left(X, d^T_N, \varepsilon \right)}{N} \leq 
     \frac{s \log L}{t}. \]
Since $s$ is an arbitrary number larger than $\dimh(X, d, \delta)$, this shows the statement.
\end{proof}

\begin{theorem} \label{theorem: mean Hausdorff dimension and entropy}
Let $(X, S, T)$ be a $\mathbb{N}^2$-action with a metric $d$ on $X$.
Suppose there exists $L>1$ such that 
\[  d\left(T x, T y\right) \leq L \, d(x,y), \quad (x,y\in X). \]
Then 
\[  \lmdimh(X, S, d) \geq \frac{\htop(X, S, T)}{\log L}. \]
\end{theorem}

\begin{proof}
For any $M>0$ we have $d^S_M\left(T x, T y\right) \leq L \cdot d^S_M(x,y)$.
Let $0<\delta<1$, $0<t<1$ and $0<\varepsilon<1$ be positive numbers with
$\delta^{1-t} < \varepsilon$.

From Lemma \ref{lemma: entropy is bounded by hausdorff dimension finitely version}, for any $M>0$
\[  \lim_{N\to \infty} \frac{\log\#\left(X, d^{S,T}_{[0,M)\times [0,N)}, \, \varepsilon\right)}{N} 
     \leq \frac{\log L}{t} \cdot \dimh\left(X, d^S_M, \delta\right).  \]
Divide the both sides by $M$ and let $M\to \infty$. We get 
\[  \lim_{\substack{N\to \infty \\ M\to \infty}} \frac{\log\#\left(X, d^{S,T}_{[0,M)\times [0,N)}, \, \varepsilon\right)}{MN} 
    \leq \frac{\log L}{t} \cdot \liminf_{M\to \infty} \frac{\dimh\left(X, d^S_M, \delta\right)}{M}. \]
Letting $\delta\to 0$, we have 
\[ \lim_{\substack{N\to \infty \\ M\to \infty}} \frac{\log\#\left(X, d^{S,T}_{[0,M)\times [0,N)}, \, \varepsilon\right)}{MN} 
    \leq \frac{\log L}{t} \cdot \lmdimh(X,S,d). \]
We can let $t\to 1$ and $\varepsilon\to 0$, so this proves the statement.
\end{proof}

\subsection{Proof of Theorem 
\ref{theorem: mean Hausdorff dimension of homogeneous systems}}
\label{subsection: proof of mean Hausdorff dimension of homogeneous systems}

Now we prove Theorem \ref{theorem: mean Hausdorff dimension of homogeneous systems}.
Recall that $\rho$ is a metric on the circle $\mathbb{R}/\mathbb{Z}$ defined by 
\[  \rho(x, y) = \min_{m\in \mathbb{Z}} |x-y-m|. \]

We write the statement of Theorem 
\ref{theorem: mean Hausdorff dimension of homogeneous systems} again.

\begin{theorem}[$=$ Theorem 
\ref{theorem: mean Hausdorff dimension of homogeneous systems}]
Let $d$ be a metric on $\left(\mathbb{R}/\mathbb{Z}\right)^\mathbb{N}$ defined by 
\[  d(x, y) = \sum_{n=1}^\infty 2^{-n}  \rho(x_n,y_n). \]
Let $a>1$ be a natural number and 
define $T_a$ on $\left(\mathbb{R}/\mathbb{Z}\right)^\mathbb{N}$ by 
the component-wise multiplication of $a$.
Let $\sigma:\left(\mathbb{R}/\mathbb{Z}\right)^\mathbb{N} \to 
\left(\mathbb{R}/\mathbb{Z}\right)^\mathbb{N}$ be the shift.
If $X\subset \left(\mathbb{R}/\mathbb{Z}\right)^\mathbb{N}$ is a closed subset satisfying 
$\sigma(X) \subset X$ and $T_a(X)\subset X$ then 
\[  \mdimh(X,\sigma,d) = \mdimm(X,\sigma,d) = \frac{\htop(X, \sigma,T_a)}{\log a}. \]
\end{theorem}

\begin{proof}
Obviously $d\left(T_a(x), T_a(y)\right) \leq a \, d(x,y)$.
So by Theorem \ref{theorem: mean Hausdorff dimension and entropy}
\[   \lmdimh(X, \sigma, d) \geq \frac{\htop(X,\sigma, T_a)}{\log a}. \]
The remaining task is to show the upper bound
\[   \umdimm(X, \sigma, d) \leq \frac{\htop(X,\sigma, T_a)}{\log a}. \]

A key fact is that, for any natural number $M$, if two points 
$u, v\in \mathbb{R}/\mathbb{Z}$ satisfy 
\[  \max_{0\leq m < M} \rho(a^m u, a^m v) < \frac{1}{2a}, \]
then 
\[   \rho(u, v) < \frac{1}{2a^M}. \]
From this, for any natural numbers $L$ and $M$, if two points $x, y\in X$ satisfy 
\[  d_{[0,L)\times[0,M)}^{\sigma,T_a}(x, y) <  \frac{1}{4a} \]
then 
\[  d(x, y) < \frac{1}{2a^M} + 2^{-L}. \]
Indeed, $d_{[0,L)\times [0,M)}^{\sigma, T_a}(x, y) <  \frac{1}{4 a}$ implies that 
for all $1\leq n \leq L$
\[ \max_{0\leq m <M} \rho(a^m x_n, a^m y_n) < \frac{1}{2a} \]
and hence 
\[  \rho(x_n, y_n) < \frac{1}{2a^M} \quad (1\leq n \leq L). \]
So 
\begin{align*}
   d(x, y) &\leq \sum_{n=1}^L 2^{-n} \rho(x_n,y_n) + \sum_{n=L+1}^\infty 2^{-n} \\
            &<  \frac{1}{2a^M} + 2^{-L}. 
\end{align*}

Let $0<\varepsilon<1$ be arbitrary. We choose natural numbers $L$ and $M$ satisfying 
\[  2^{-L} < \frac{\varepsilon}{2}, \quad   a^{-M}  \leq \varepsilon < a^{-M+1}. \]
From the above consideration, for any natural number $N$
\[   d^{\sigma, T_a}_{[0, N+L)\times [0,M)}(x,y) < \frac{1}{4a} \Longrightarrow 
    d_N^{\sigma}(x,y) < \frac{1}{2a^M} + 2^{-L} < \varepsilon. \]
So
\[  \#\left(X, d^\sigma_N,\varepsilon\right) \leq 
    \#\left(X, d^{\sigma, T_a}_{[0, N+L)\times [0,M)},\, \frac{1}{4a}\right). \]
Hence 
\begin{align*}
    \lim_{N\to \infty} \frac{\log \#\left(X, d^\sigma_N,\varepsilon\right)}{N}  &\leq 
    \lim_{N\to \infty} 
    \frac{\log \#\left(X, d^{\sigma, T_a}_{[0, N+L)\times [0,M)},\, \frac{1}{4a}\right)}{N} \\
   & = \lim_{N\to \infty} 
    \frac{\log \#\left(X, d^{\sigma, T_a}_{[0, N)\times [0,M)},\, \frac{1}{4a}\right)}{N}.
\end{align*}

From $\varepsilon < a^{-M+1}$, we have $(M-1)\log a < \log(1/\varepsilon)$ and 
\[  \frac{1}{\log(1/\varepsilon)} < \frac{1}{(M-1)\log a} =
     \frac{M}{(M-1)\log a}\cdot \frac{1}{M}. \]
Therefore 
\begin{align*}
  \lim_{N\to \infty} \frac{\dimm\left(X, d^\sigma_N,\varepsilon\right)}{N} & =
  \lim_{N\to \infty} 
  \frac{\log \#\left(X, d^\sigma_N,\varepsilon\right)}{N \log (1/\varepsilon)}  \\
  & \leq  \frac{M}{(M-1)\log a}\cdot
    \lim_{N\to \infty} 
    \frac{\log \#\left(X, d^{\sigma, T_a}_{[0, N)\times [0,M)},\, \frac{1}{4a}\right)}{NM}. 
\end{align*}    
$M$ goes to infinity as $\varepsilon$ goes to zero.
So 
\begin{align*}
   \limsup_{\varepsilon\to 0} \left(\lim_{N\to \infty} 
   \frac{\dimm\left(X, d^\sigma_N,\varepsilon\right)}{N}\right)
   & \leq  \frac{1}{\log a} \cdot \lim_{\substack{N\to \infty\\ M\to \infty}} 
    \frac{\log \#\left(X, d^{\sigma, T_a}_{[0, N)\times [0,M)},\, \frac{1}{4a}\right)}{NM} \\
   & \leq \frac{1}{\log a} \cdot \htop(X, \sigma, T_a).
\end{align*}
This proves $\umdimm(X,\sigma,d) \leq \frac{\htop(X,\sigma,T_a)}{\log a}$.
\end{proof}

\section{Self-similarity and mean Hausdorff dimension}  \label{section: self-similarity and mean Hausdorff dimension}

In this section we introduce a “self-similar system”, which is a dynamical version of a self-similar set.
We show that mean Hausdorff dimension and metric mean dimension coincide for such systems.

\subsection{Self-similar systems}  \label{subsection: self-similar systems}

Let 
\[ \ell^\infty = \left\{(x_n)_{n\in \mathbb{N}}\in \mathbb{R}^{\mathbb{N}}\, \middle|\, \sup_{n\geq 1} |x_n| < \infty\right\} \]
be the space of bounded sequences 
with the norm $\norm{x}_\infty := \sup_{n\geq 1} |x_n|$ for 
$x = (x_n)_{n\in \mathbb{N}}$.
We always assume that $\ell^\infty$ is endowed with the weak$^*$ topology as the dual space of $\ell^1$.
We define the shift map $\sigma:\ell^\infty \to \ell^\infty$ by
\[  \sigma\left((x_n)_{n\in \mathbb{N}}\right) = (x_{n+1})_{n\in \mathbb{N}}. \]
This is continuous with respect to the weak$^*$ topology.
We will consider a certain self-similar set of $\ell^\infty$ invariant under the shift map $\sigma$.

Let $(\Omega, T)$ be a dynamical system.
Suppose that for each $\omega\in \Omega$ we are given a point 
$a(\omega) = \left(a(\omega)_n\right)_{n\in \mathbb{N}}\in \ell^\infty$ so that the map 
\[  \Omega \ni  \omega \mapsto a(\omega) \in \ell^\infty \]
is continuous (with respect to the weak$^*$ topology of $\ell^\infty$) and equivariant (i.e. $\sigma\left(a(\omega)\right)  = a(T\omega)$).
Since $\Omega$ is compact, we have 
\[   \sup_{\omega \in \Omega} \norm{a(\omega)}_\infty < \infty. \]
Fix a real number $c$ with $0<c<1$.
For each $\omega \in \Omega$ we define a contracting similarity transformation $S_\omega:\ell^\infty\to \ell^\infty$ by
\[  S_\omega (x) = c x + a(\omega). \]
Then $\sigma\left(S_\omega(x)\right) = S_{T\omega}\left(\sigma(x)\right)$.

\begin{proposition}[Definition of a self-similar system] \label{prop: definition of self-similar system}
There uniquely exists a non-empty compact subset $X$ of $\ell^\infty$ satisfying 
\[   X = \bigcup_{\omega \in \Omega} S_\omega(X). \]
This $X$ is $\sigma$-invariant (i.e. $\sigma(X) \subset X$).
The dynamical system $(X, \sigma)$ is called a \textbf{self-similar system} defined by the family of 
contracting similarity transformations $\{S_\omega\}_{\omega\in \Omega}$.
\end{proposition}

Notice that here $X$ is compact with respect to the weak$^*$ topology, not the norm topology,
and that $X$ becomes bounded and closed in the $\ell^\infty$-norm
by the uniform boundedness principle (the Banach--Steinhaus theorem).

\begin{proof}
We define $X\subset \ell^\infty$ by 
\[  X = \left\{\sum_{k=0}^\infty c^k a(\omega_k) \, \middle|\, \omega_k\in \Omega \> (k\geq 0)\right\}. \] 
Notice that the sum $\sum_{k=0}^\infty c^k a(\omega_k)$ absolutely converges because 
$\sup_{\omega\in \Omega}\norm{a(\omega)}_\infty < \infty$.

$X$ is the image of a continuous map 
\[  \Omega\times \Omega\times \Omega\times \cdots \to \ell^\infty, \quad 
     (\omega_k)_{k\geq 0} \mapsto \sum_{k=0}^\infty c^k a(\omega_k), \]
where $\Omega\times \Omega\times \cdots$ is endowed with the product topology.
Since $\Omega$ is compact, $X$ is also compact.
We have 
\[   S_\omega\left(\sum_{k=0}^\infty c^k a(\omega_k)\right) = a(\omega) + \sum_{k=0}^\infty c^{k+1}a(\omega_k). \]
From this, it is easy to see that $X$ satisfies 
\[  X =  \bigcup_{\omega \in \Omega} S_\omega(X). \]
We also have 
\[  \sigma\left(\sum_{k=0}^\infty c^k a(\omega_k)\right)  = \sum_{k=0}^\infty c^k \sigma\left(a(\omega_k)\right) 
     = \sum_{k=0}^\infty c^k a\left(T\omega_k\right). \]
Therefore $\sigma(X) \subset X$.

Next we study the uniqueness of $X$.
Suppose a non-empty compact subset $Y\subset \ell^\infty$ satisfies 
\[  Y = \bigcup_{\omega\in \Omega} S_\omega(Y). \]
Recall that the compactness (with respect to the weak$^*$ topology) 
implies that $Y$ is bounded and closed in the $\ell^\infty$-norm
by the uniform boundedness principle.

Suppose $Y\not \subset X$. We set 
\[  \delta := \sup_{y \in Y}\left(\inf_{x\in X} \norm{x-y}_\infty \right) >0. \]
Since $Y = \bigcup_{\omega\in \Omega} S_\omega(Y)$, we have 
\[  \delta = \sup_{\substack{y\in Y\\ \omega\in \Omega}}\left(\inf_{x\in X} \norm{x-S_\omega y}_\infty\right). \]
However 
\begin{align*}
     \inf_{x\in X} \norm{x-S_\omega y}_\infty & \leq \inf_{x\in X} \norm{S_\omega x - S_\omega y}_\infty 
     \quad (\text{by $S_\omega(X) \subset X$}) \\
      &= c\cdot \inf_{x \in X} \norm{x-y}_\infty \leq c\cdot \delta. 
\end{align*}      
Hence $\delta\leq c\cdot \delta < \delta$. This is a contradiction. Hence $Y\subset X$.
By switching the roles of $X$ and $Y$, we also have $X \subset Y$.
So $X=Y$. This shows the uniqueness of $X$.
\end{proof}

\begin{remark}
We recall that the map $S_\omega:\ell^\infty \to \ell^\infty$ has the form
\[ S_\omega(u)  = c u + a(\omega), \quad (\text{$c$ is a fixed constant}). \]
The point is that the linear part of $S_\omega$ is a scalar multiplication.
Probably many readers feel that this form is too restricted.
This severe restriction comes from the fact that the $\ell^\infty$-geometry does not admit many similarity transformations.
There is no “rotation” for $\ell^\infty$ (except for permutations of coordinates).
We expect that it is more interesting 
to study “self-affine sets” of $\ell^\infty$, rather than “self-similar sets”,
because the space $\ell^\infty$ seems to admit many interesting affine transformations.
We hope to come back to this study in a future paper.
\end{remark}

Let $X\subset \ell^\infty$ be the self-similar system defined by $\{S_\omega\}_{\omega\in \Omega}$.
We define a metric $d$ on $X$ by 
\[  d\left((x_n)_{n\in \mathbb{N}}, (y_n)_{n\in \mathbb{N}}\right) = \sum_{n=1}^\infty 2^{-n} |x_n-y_n|. \]
We are interested in the mean Hausdorff dimension and metric mean dimension of 
$(X, \sigma, d)$.
The following is the main result of this section.

\begin{theorem} \label{main theorem for self-similar systems}
Under the above setting, the self-similar system $X$ satisfies 
\[  \mdimh(X,\sigma,d) = \mdimm(X, \sigma, d)   \leq \frac{\htop(\Omega, T)}{\log (1/c)}. \]
Here $\htop(\Omega, T)$ is the topological entropy of $(\Omega, T)$.
\end{theorem}

Therefore, the mean Hausdorff dimension and metric mean dimension coincide for self-similar systems.

\begin{remark}
Let $a_1, \dots, a_m$ be vectors in $\mathbb{R}^N$, and define contracting similarity transformations
$f_i:\mathbb{R}^N\to \mathbb{R}^N$ $(1\leq i \leq m)$ by 
\[  f_i(u) = cu + a_i. \]
Let $K\subset \mathbb{R}^N$ be an attractor of $\{f_1, \dots, f_m\}$.
Then the \textit{similarity dimension} of $K$ is given by 
\[  \frac{\log m}{\log(1/c)}. \]
It is well-known that the Hausdorff and Minkowski dimensions of $K$ (with respect to the Euclidean metric) are 
bounded by the similarity dimension.
The term $\frac{\htop(\Omega, T)}{\log (1/c)}$ in Theorem \ref{main theorem for self-similar systems} is an analogue of 
the similarity dimension in our setting.
\end{remark}

We prove Theorem \ref{main theorem for self-similar systems} in the next subsection.
Before going into the proof, we study a simple example.

\begin{example}[$\beta$-expansions] \label{example: beta-expansions}
Let $a>1$ be an integer greater than one.
Let $\{0, 1, 2, \dots, a-1\}^{\mathbb{N}}$ be the full-shift on the alphabet $0,1,2,\dots, a-1$.
We naturally consider that this is a subset of $\ell^\infty$:
\[  \{0,1,2,\dots, a-1\}^{\mathbb{N}}\subset \ell^\infty. \]
Let $\Omega\subset \{0,1,2,\dots, a-1\}^{\mathbb{N}}$ be a subshift 
(a closed subset invariant under the shift map $\sigma$).
We fix a real number $\beta$ with $\beta\geq a$.
For each $\omega \in \Omega$ we define $S_\omega:\ell^\infty \to \ell^\infty$ by 
\[   S_\omega(x) = \frac{x+\omega}{\beta}. \]
Let $X\subset \ell^\infty$ be a self-similar system defined by these transformations.
This is given by\footnote{Notice that, when $\beta = a$, 
the system $(X,\sigma)$ provides an example for 
Theorem \ref{theorem: mean Hausdorff dimension of homogeneous systems} by projecting it to 
the infinite dimensional torus $\left(\mathbb{R}/\mathbb{Z}\right)^{\mathbb{N}}$. 
We also note that, when $\beta>a$, 
our setting is simpler than general $\beta$-expansions because we restrict “digits” to $\{0,1,2,\dots, a-1\}$ and
$a$ is forbidden.
For example, when $\beta = \frac{1+\sqrt{5}}{2} = 1.618\dots$, a complication of $\beta$-expansions may occur from 
$\frac{1}{\beta} = \frac{1}{\beta^2} + \frac{1}{\beta^3}$.
We sidestep this complication simply by forbidding the digit $1$ to appear.
We do not dig deeper into this problem in this paper.}
\[  X = \left\{\sum_{n=1}^\infty \frac{\omega_n}{\beta^n}\middle|\, \omega_n\in \Omega\right\}. \]
By Theorem \ref{main theorem for self-similar systems} we have 
\[   \mdimh(X,\sigma,d) = \mdimm(X, \sigma, d)   \leq \frac{\htop(\Omega, \sigma)}{\log \beta}. \]
Actually the equality holds here as we will see below.
We need the next claim.

\begin{claim} \label{claim: beta expansion}
Let $0\leq u_k, v_k\leq a-1$ be integers $(1\leq k \leq n)$.
If $(u_1, \dots, u_n) \neq (v_1, \dots, v_n)$ then 
\[    \left|\sum_{k=1}^n \frac{u_k}{\beta^k} - \sum_{k=1}^n \frac{v_k}{\beta^k}\right| \geq \frac{1}{\beta^n}. \]
\end{claim}

\begin{proof}
Take an integer $m \in [1,n]$ with 
$(u_1, \dots, v_{m-1}) = (v_1, \dots, v_{m-1})$ and $u_m\neq v_m$.
We have 
\begin{align*}
     \left|\sum_{k=1}^n \frac{u_k}{\beta^k} - \sum_{k=1}^n \frac{v_k}{\beta^k}\right|  &=
      \left|\frac{u_m-v_m}{\beta^m} + \sum_{k=m+1}^n \frac{u_k-v_k}{\beta^k}\right| \\
      &\geq \frac{|u_m-v_m|}{\beta^m} - \sum_{k=m+1}^n \frac{|u_k-v_k|}{\beta^k} \\
      &\geq \frac{1}{\beta^m} - \sum_{k=m+1}^n \frac{a-1}{\beta^k} \\
      &= \frac{1}{\beta^m} - (a-1) \cdot \frac{\beta^{n-m}-1}{\beta^n(\beta-1)} \\
      & \geq \frac{1}{\beta^m} - \frac{\beta^{n-m}-1}{\beta^n} \quad 
      (\text{by $a\leq \beta$}) \\
     & = \frac{1}{\beta^n}.
\end{align*}
\end{proof}
Let $N$ be a natural number and let $\pi_N:\{0,1,2,\dots, a-1\}^{\mathbb{N}}\to \{0,1,2,\dots, a-1\}^N$ be the projection 
to the first $N$ coordinates.
Fix $\xi \in \Omega$. 
For $\omega_1, \dots, \omega_n, \omega^\prime_1, \dots, \omega^\prime_n \in \Omega$, if 
$\left(\pi_N(\omega_1), \dots, \pi_N(\omega_n)\right) \neq \left(\pi_N(\omega^\prime_1), \dots, \pi_N(\omega^\prime_n)\right)$
then by Claim \ref{claim: beta expansion}
\[  d_N\left(\sum_{k=1}^n \frac{\omega_k}{\beta^k} + \sum_{k=n+1}^\infty \frac{\xi}{\beta^k}, 
                 \sum_{k=1}^n \frac{\omega^\prime_k}{\beta^k} + \sum_{k=n+1}^\infty \frac{\xi}{\beta^k}\right) \geq \frac{1}{\beta^n}. \]
Therefore 
\[  \#\left(X, d_N, \frac{1}{\beta^n}\right) \geq \left|\pi_N(\Omega)\right|^n.  \]
Then 
\[  \frac{\dimm\left(X, d_N, \frac{1}{\beta^n}\right)}{N} \geq \frac{\log \left|\pi_N(\Omega)\right|}{N\log \beta}. \]
Letting $N\to \infty$, the right-hand side becomes $\htop(\Omega, \sigma)/\log \beta$:
\[  \lim_{N\to \infty} \frac{\dimm\left(X, d_N, \frac{1}{\beta^n}\right)}{N} \geq 
           \frac{\htop(\Omega, \sigma)}{\log \beta}. \]
Letting $n\to \infty$, we conclude 
\[   \mdimm(X, \sigma, d) \geq    \frac{\htop(\Omega, \sigma)}{\log \beta}. \]
Thus we have 
\[ \mdimh(X,\sigma,d) = \mdimm(X, \sigma, d)   =\frac{\htop(\Omega, \sigma)}{\log \beta}. \]
\end{example}

\subsection{Proof of Theorem \ref{main theorem for self-similar systems}}
\label{subsection: proof of main theorem for self-similar systems}

The following is a key lemma for the proof of Theorem \ref{main theorem for self-similar systems}.
This is a finite accuracy version of a theorem of Falconer \cite[Theorem 4]{Falconer}.
The proof closely follows Falconer’s original argument.

\begin{lemma} \label{lemma: finite accuracy Falconer lemma}
For any real numbers $0<\varepsilon, a, \tau<1$ there exists 
$\delta_0 = \delta_0(\varepsilon,a,\tau)>0$ such that the following statement holds true.
Let $(X, d)$ be a compact metric space such that for every closed ball $B \subset X$ of radius $\varepsilon$
there exists a map $\varphi:X\to B$ satisfying 
\[  d\left(\varphi(x), \varphi(y)\right) \geq a\varepsilon d(x, y) \quad (x, y\in X). \]
Then
\[   \dimh(X, d, \delta_0) \geq \tau \cdot \frac{\log \#(X, d, 9\varepsilon)}{\log\left(\frac{1}{a\varepsilon}\right)}. \]
\end{lemma}

\begin{proof}
Set $N = \#(X, d, 9\varepsilon)$. There exists points $x_1, \dots, x_N\in X$ such that 
$d(x_i, x_j)> 3\varepsilon$ for $i \neq  j$.
Let $B_i$ be the closed ball of radius $\varepsilon$ centered at $x_i$. Then 
\[  d\left(B_i, B_j\right) := \inf_{\substack{x\in B_i \\ y\in B_j}} d(x, y) > \varepsilon \quad (>a\varepsilon). \]
We can take a map $\varphi_i:X\to B_i$ satisfying 
$d\left(\varphi_i(x), \varphi_i(y)\right) \geq a \varepsilon d(x, y)$.

For $1\leq i_1, i_2, \dots, i_n\leq N$ we set 
\[   B_{i_1 i_2 \dots i_n} = \varphi_{i_n}\circ \varphi_{i_{n-1}}\circ \cdots \circ \varphi_{i_2}(B_{i_1}). \]
We have 
\[  B_{i_1 i_2 \dots i_n}\subset B_{i_2 i_3 \dots i_n} \subset 
     \dots \subset B_{i_{n-1} i_n} \subset B_{i_n}. \]
See Figure \ref{figure: nested balls}.

\begin{figure}[htbp]
\begin{center}
\includegraphics[scale=0.7]{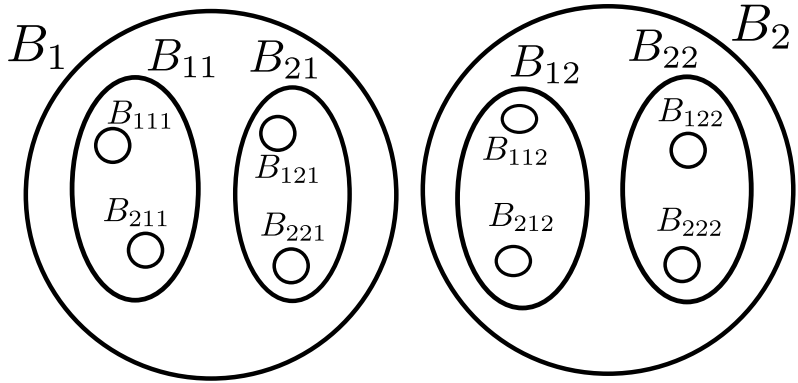} 
\caption{$B_{i_1 i_2 i_3}\subset B_{i_2 i_3} \subset B_{i_3}$ for $1\leq i_1, i_2, i_3\leq 2$.}  \label{figure: nested balls}
\end{center}
\end{figure}

\begin{claim}   \label{claim: distance between nested balls}
For $(i_1, i_2, \dots, i_n) \neq (j_1, j_2, \dots, j_n)$,
\[  d\left(B_{i_1 i_2 \dots i_n}, B_{j_1 j_2 \dots j_n}\right) > (a\varepsilon)^n. \]
\end{claim}
\begin{proof}
If $i_n\neq j_n$ then $B_{i_1 \dots i_n} \subset B_{i_n}$ and $B_{j_1 \dots j_n}\subset B_{j_n}$ imply
\[  d\left(B_{i_1 \dots i_n}, B_{j_1 \dots j_n}\right) \geq  d(B_{i_n}, B_{j_n}) > a\varepsilon \geq  (a\varepsilon)^n. \]
If $i_n = j_n$ and $i_{n-1}\neq j_{n-1}$ then 
$B_{i_1 \dots i_n} \subset \varphi_{i_n}(B_{i_{n-1}})$ and 
$B_{j_1 \dots j_n}\subset \varphi_{i_n}(B_{j_{n-1}})$ imply 
\begin{align*}
  d\left(B_{i_1 \dots i_n}, B_{j_1 \dots j_n}\right) & \geq d\left( \varphi_{i_n}(B_{i_{n-1}}), \varphi_{i_n}(B_{j_{n-1}})\right) \\
    &\geq  a\varepsilon d\left(B_{i_{n-1}}, B_{j_{n-1}}\right) \\
    &  > (a\varepsilon)^2 \geq (a\varepsilon)^n. 
\end{align*}    
If $(i_n, i_{n-1}) = (j_n, j_{n-1})$ and $i_{n-2}\neq j_{n-2}$ then 
$B_{i_1 \dots i_n} \subset \varphi_{i_n}\circ \varphi_{i_{n-1}}(B_{i_{n-2}})$ and 
$B_{j_1 \dots j_n}\subset \varphi_{i_n}\circ \varphi_{i_{n-1}}(B_{j_{n-2}})$ imply 
\begin{align*}
     d\left(B_{i_1 \dots i_n}, B_{j_1 \dots j_n}\right) & \geq  
     d\left(\varphi_{i_n}\circ \varphi_{i_{n-1}}(B_{i_{n-2}}), \varphi_{i_n}\circ \varphi_{i_{n-1}}(B_{j_{n-2}})\right) \\
     & \geq a\varepsilon d\left(\varphi_{i_{n-1}}(B_{i_{n-2}}), \varphi_{i_{n-1}}(B_{j_{n-2}})\right) \\
     & \geq (a\varepsilon)^2 d\left(B_{i_{n-2}}, B_{j_{n-2}}\right) \\
     & \geq (a\varepsilon)^3 \geq (a\varepsilon)^n.
\end{align*}
We can proceed similarly and prove the claim.
\end{proof}

For each $n\geq 1$ we take a Borel probability measure $\mu_n$ on $X$ such that for every 
$1\leq i_1, \dots, i_n\leq N$ we have 
\[  \mu_n\left(B_{i_1 \dots i_n}\right) = \frac{1}{N^n}. \]
Notice that this implies 
\[  \mu_n\left(\bigcup_{1\leq i_1, \dots, i_n\leq N} B_{i_1 \dots i_n}\right) = 1. \]
Moreover, for any $1\leq m \leq n$ and $1\leq i_1, \dots, i_m \leq N$ we also have 
\[ \mu_n\left(B_{i_1 \dots i_m}\right) = \frac{1}{N^m}. \]

We fix $0<\delta_0<1$ satisfying 
\begin{equation} \label{eq: choice of delta_0}
     (1-\tau)\cdot \frac{\log \delta_0}{\log(a\varepsilon)} \geq 1.
\end{equation}     
Set 
\[   s = \frac{\tau \log N}{\log\left(\frac{1}{a\varepsilon}\right)}. \]

\begin{claim} \label{claim: power law for mu_n}
Let $E\subset X$ be a Borel subset satisfying $0< \diam \,E < \delta_0$ then 
\[  \mu_n(E) \leq \left(\diam\, E\right)^{s} \]
for all sufficiently large $n$.
\end{claim}

\begin{proof}
Set $r = \diam\, E$.
Take a natural number $m$ with $(a\varepsilon)^{m+1} \leq r < (a\varepsilon)^m$.
Then $E$ intersects with at most one set in $\left\{B_{i_1 i_2 \dots i_m}\right\}$
by Corollary \ref{claim: distance between nested balls}.
For any integer $n\geq m$ we have 
\[  \mu_n(E) \leq \frac{1}{N^m}. \]
From $(a\varepsilon)^{m+1}\leq r$, we have $(m+1)\log (a\varepsilon) \leq \log r$ and hence 
\[  m \geq -1 + \frac{\log r}{\log(a\varepsilon)}. \]
From (\ref{eq: choice of delta_0})
\begin{align*}
   m &\geq -(1-\tau)\cdot  \frac{\log \delta_0}{\log(a\varepsilon)}  + \frac{\log r}{\log(a\varepsilon)} \\
       &> -(1-\tau)\cdot  \frac{\log r}{\log(a\varepsilon)}  + \frac{\log r}{\log(a\varepsilon)} \quad (\text{by $r<\delta_0$}) \\
       &= \tau\cdot \frac{\log r}{\log(a\varepsilon)}.
\end{align*}
Then 
\begin{align*}
  \frac{1}{N^m} & = \exp\left(-m\log N\right) \\
   & < \exp\left(-\tau\cdot \frac{\log r}{\log(a\varepsilon)}\cdot \log N\right) \\
   & = \exp\left(s \log r\right)  \quad 
   \left(\text{by } s = \frac{\tau \log N}{\log\left(\frac{1}{a\varepsilon}\right)}\right) \\
   & = r^s.
\end{align*}
Thus we have $\mu_n(E) < r^s$ for any $n\geq m$.
\end{proof}

Suppose we are given an open cover $X = U_1\cup \dots \cup U_K$ with 
$0<\diam\, U_k < \delta_0$ for all $1\leq k \leq K$.
Then by Claim \ref{claim: power law for mu_n}
\[  \mu_n(U_k) < \left(\diam \, U_k\right)^s \]
for all $k$ and any sufficiently large $n$. Therefore 
\[  \sum_{k=1}^K \left(\diam\, U_k\right)^s > \sum_{k=1}^K \mu_n(U_k) \geq \mu_n(X) = 1. \]
This implies $\mathcal{H}^s_{\delta_0}(X, d) \geq 1$.
Thus 
\[  \dimh(X, d, \delta_0) \geq s = \tau \cdot \frac{\log \#(X, d, 9\varepsilon)}{\log\left(\frac{1}{a\varepsilon}\right)}. \]
This proves the statement of the lemma.
\end{proof}

We return to the setting of \S \ref{subsection: self-similar systems}.
First we recall various definitions.
Let $(\Omega, T)$ be a dynamical system.
Suppose that we are given a continuous equivariant map 
\[  \Omega \ni \omega \mapsto a(\omega)\in \ell^\infty. \] 
Fix $0<c<1$ and define $S_\omega:\ell^\infty \to \ell^\infty$ for each $\omega\in \Omega$ by 
\[  S_\omega(x) = c x + a(\omega). \]

Let $X\subset \ell^\infty$ be the self-similar system defined by $\{S_\omega\}_{\omega\in \Omega}$.
It is a shift-invariant, non-empty, compact subset of $\ell^\infty$ 
(with respect to the weak$^*$ topology) satisfying 
\[   X= \bigcup_{\omega\in \Omega} S_\omega(X). \]
We define a metric $d$ on $X$ by 
\[  d\left((x_n)_{n\in \mathbb{N}}, (y_n)_{n\in \mathbb{N}}\right) = \sum_{n=1}^\infty 2^{-n} |x_n-y_n|. \]
We are interested in the mean Hausdorff dimension and metric mean dimension of $(X, \sigma, d)$,
where $\sigma:X\to X$ is the shift map.

Since $X$ is compact, we can find $A\geq 1$ such that
\[   \diam(X, d) < A. \]
For $N\geq 1$ we define a metric $d_N$ on $X$ by 
\[   d_N(x, y) = \max_{0\leq n < N} d\left(\sigma^n(x), \sigma^n(y)\right). \]

\begin{claim} \label{claim: self-similar copy in balls}
Let $0<\varepsilon<1$ and let $N$ be a natural number.
For any point $p\in X$ there exists a map $\varphi:X\to B_\varepsilon(p, d_N)$ satisfying 
\[  d_N\left(\varphi(x), \varphi(y)\right) \geq \frac{c \varepsilon}{A} d_N(x,y). \]
Here $B_\varepsilon(p, d_N)$ is the closed $\varepsilon$-ball centered at $p$ with respect to the metric $d_N$.
\end{claim}

\begin{proof}
For any $\omega\in \Omega$
\[
  d\left(S_\omega(x), S_\omega(y)\right) = d\left(cx + a(\omega), cy+a(\omega)\right) = \sum_{n=1}^\infty 2^{-n}|c x_n - c y_n| 
   = c d(x, y).
\]
Similarly we have
\begin{equation} \label{eq: similarity condition for d_N}
   d_N\left(S_\omega(x), S_\omega(y)\right) = c d_N(x,y). 
\end{equation}   

Take a natural number $k$ with $c^k A\leq \varepsilon <  c^{k-1}A$.
From $X= \bigcup_{\omega\in \Omega} S_\omega(X)$ we can find a sequence $\omega_1, \dots, \omega_k\in \Omega$ satisfying 
\[  p \in S_{\omega_1}\circ \dots \circ S_{\omega_k}(X). \]
Set $\varphi := S_{\omega_1}\circ \dots \circ S_{\omega_k}: X\to X$.
For any $x,y\in X$, by (\ref{eq: similarity condition for d_N})
\[  d_N\left(\varphi(x), \varphi(y)\right) = c^k d_N(x,y) \leq c^k A \leq \varepsilon. \]
Therefore $\varphi(X)\subset B_\varepsilon(p, d_N)$.
From $\varepsilon <  c^{k-1}A$ we have $c^k > \frac{c \varepsilon}{A}$ and hence 
\[  d_N\left(\varphi(x), \varphi(y)\right) = c^k d_N(x,y) \geq \frac{c\varepsilon}{A} d_N(x, y). \]
\end{proof}

Now we start the proof of Theorem \ref{main theorem for self-similar systems}.
We rewrite the statement:

\begin{theorem}[$=$ Theorem \ref{main theorem for self-similar systems}]
For the above self-similar system $(X, \sigma)$ we have
\[   \mdimh(X, \sigma, d)  = \mdimm(X,\sigma, d) \leq \frac{\htop(\Omega, T)}{\log (1/c)}. \]
\end{theorem}

\begin{proof}
We first prove $\mdimh(X, \sigma, d)  = \mdimm(X,\sigma, d)$.
It is enough to show 
\[  \lmdimh(X, \sigma, d) \geq  \umdimm(X, \sigma, d). \]
Let $0<\varepsilon, \tau<1$.
By Lemma \ref{lemma: finite accuracy Falconer lemma} and Claim \ref{claim: self-similar copy in balls}, there exists 
$\delta_0 = \delta_0(\varepsilon, c/A, \tau)>0$ such that for any natural number $N$
\[  \dimh(X, d_N, \delta_0) \geq \tau\cdot \frac{\log \#(X, d_N, 9\varepsilon)}{\log\left(\frac{A}{c\varepsilon}\right)}. \]
Divide this by $N$ and let $N\to \infty$:
\[ \lmdimh(X, \sigma, d) \geq  \liminf_{N\to \infty}  \frac{\dimh(X, d_N, \delta_0)}{N} \geq  \lim_{N\to \infty}
    \tau \cdot \frac{\log \#(X, d_N, 9\varepsilon)}{N \log\left(\frac{A}{c\varepsilon}\right)}.  \]
Letting $\varepsilon\to 0$, we get 
\[  \lmdimh(X, \sigma, d)  \geq \tau \cdot \umdimm(X, \sigma, d). \]
Letting $\tau\to 1$, we conclude:
\[  \lmdimh(X, \sigma, d) \geq  \umdimm(X, \sigma, d). \]

Next we prove 
\[  \mdimm(X,\sigma, d) \leq \frac{\htop(\Omega, T)}{\log (1/c)}. \]
Let $\rho$ be a metric on $\Omega$.
For $N\geq 1$ we define a metric $\rho_N$ on $\Omega$ by 
\[  \rho_N(\omega, \omega^\prime) = \max_{0\leq n < N} \rho\left(T^n \omega, T^n \omega^\prime\right). \]

For $\omega, \omega^\prime\in \Omega$ and $x, y\in X$
\begin{align*}
   d\left(S_\omega(x), S_{\omega^\prime}(y)\right)  & = d\left(c x + a(\omega), c y+ a(\omega^\prime)\right) \\
      & = \sum_{n=1}^\infty 2^{-n} \left|c x_n + a(\omega)_n - c y_n - a(\omega^\prime)_n\right| \\
      & \leq c d(x, y) + d\left(a(\omega), a(\omega^\prime)\right).
\end{align*}
Similarly, for any natural number $N$
\begin{equation*}
  d_N\left(S_\omega(x), S_{\omega^\prime}(y)\right)  \leq c d_N(x, y) + d_N\left(a(\omega), a(\omega^\prime)\right). 
\end{equation*}  
By repeatedly applying this inequality,
for $\omega_1, \dots, \omega_n, \omega^\prime_1, \dots, \omega^\prime_n \in \Omega$ and $x, y\in X$
\begin{equation}   \label{eq: distortion bound}
   \begin{split}
   d_N\left(S_{\omega_1}\circ \dots \circ S_{\omega_n}(x), S_{\omega^\prime_1}\circ \dots \circ S_{\omega^\prime_n}(y)\right) &\leq 
   c^n d_N(x, y) + \sum_{i=1}^n c^{i-1} d_N\left(a(\omega_i), a(\omega^\prime_i)\right) \\
  &\leq  c^n d_N(x, y) + \frac{\max_{1\leq i \leq n} d_N\left(a(\omega_i), a(\omega^\prime_i)\right)}{1-c}  \\
   & < c^n A +    \frac{\max_{1\leq i \leq n} d_N\left(a(\omega_i), a(\omega^\prime_i)\right)}{1-c}.
   \end{split}
\end{equation}
In the last inequality we have used $\diam(X, d_N) = \diam(X, d) < A$.

Let $0<\varepsilon<1$. We take $\delta>0$ such that 
\[  \rho(\omega, \omega^\prime) < \delta \Longrightarrow  d\left(a(\omega), a(\omega^\prime)\right) < (1-c) \frac{\varepsilon}{6}. \]
Then for any natural number $N$ we have 
\[  \rho_N(\omega, \omega^\prime) < \delta \Longrightarrow  d_N\left(a(\omega), a(\omega^\prime)\right) < (1-c) \frac{\varepsilon}{6}. \]
We take a natural number $n$ satisfying 
\[   Ac^n < \frac{\varepsilon}{6} \leq Ac^{n-1}. \]

Let $\Omega_{N,\delta}\subset \Omega$ be a $\delta$-spanning set with respect to $\rho_N$ with 
$\left|\Omega_{N,\delta}\right|  = \#\left(\Omega, \rho_N, \delta\right)$.
Fix $p\in X$.
For any $\omega_1, \dots, \omega_n\in \Omega$ we can find $\omega^\prime_1, \dots, \omega^\prime_n\in \Omega_{N,\delta}$
satisfying $\rho_{N}(\omega_i, \omega^\prime_i) < \delta$ (and then 
$d_N\left(a(\omega_i), a(\omega^\prime_i)\right) < (1-c) \varepsilon/6$).
Then for any $x\in X$, by (\ref{eq: distortion bound})
\[  d_N\left(S_{\omega_1}\circ \dots \circ S_{\omega_n}(x), S_{\omega^\prime_1}\circ \dots \circ S_{\omega^\prime_n}(p)\right)
     < \frac{\varepsilon}{3}. \]
Since we know that 
\[  X = \bigcup_{\omega_1, \dots, \omega_n\in \Omega} S_{\omega_1}\circ \dots \circ S_{\omega_n}(X), \]
this implies that the set 
\[  \left\{S_{\omega^\prime_1} \circ\dots \circ S_{\omega^\prime_n}(p) \middle|\, 
      \omega^\prime_1, \dots, \omega^\prime_n\in \Omega_{N,\delta}\right\} \]
is a $\varepsilon/3$-spanning set of $X$ with respect to $d_N$.
Therefore
\[  \#\left(X, d_N, \varepsilon\right) \leq \left|\Omega_{N,\delta}\right|^n = \left\{\#\left(\Omega, \rho_N, \delta\right)\right\}^n. \]
Since $\varepsilon/6 \leq A c^{n-1}$, we have 
\[  n\log(1/c) \leq \log\left(\frac{6A}{c\varepsilon}\right). \]
Hence 
\[  \frac{\log\#\left(X, d_N,\varepsilon\right)}{N \log\left(\frac{6A}{c\varepsilon}\right)} 
     \leq  \frac{\log \#\left(\Omega, \rho_N, \delta\right)}{N \log (1/c)}. \]
Letting $N\to \infty$ 
\[  \lim_{N\to \infty}  \frac{\log\#\left(X, d_N,\varepsilon\right)}{N \log\left(\frac{6A}{c\varepsilon}\right)} 
    \leq \lim_{N\to \infty} \frac{\log \#\left(\Omega, \rho_N, \delta\right)}{N \log (1/c)}  \leq 
          \frac{\htop(\Omega, T)}{\log(1/c)}. \]
Letting $\varepsilon\to 0$, we conclude 
\[  \mdimm(X, \sigma, d) \leq  \frac{\htop(\Omega, T)}{\log(1/c)}. \]
\end{proof}

\section{Infinite dimensional carpets}  \label{section: infinite dimensional carpets}

In this section we study a mean dimension version of Bedford--McMullen carpets.
This provides a natural example for which mean Hausdorff dimension and metric mean dimension do not 
coincide.

\subsection{Weighted topological entropy} \label{subsection: weighted topological entropy}

Here we review the theory of weighted topological entropy.
We will need this notion for formulating the main result of this section.
The weighted topological entropy was originally introduced by Feng--Huang \cite{Feng--Huang}.
The presentation here follows the approach of \cite{Tsukamoto_weighted}.

Let $(X, T)$ and $(Y, S)$ be dynamical systems.
A map $\pi:X\to Y$ is called a \textbf{factor map} if $\pi$ is a continuous surjection satisfying 
$\pi\circ T = S\circ \pi$.
We often use the notation $\pi:(X, T)\to (Y, S)$ for clarifying the underlying dynamics.

For a real number $0\leq w\leq 1$ and a factor map $\pi:(X, T) \to (Y, S)$ we will define a 
weighted topological entropy $\htop^w(\pi,X, T)$.
Let $d$ and $d^\prime$ be metrics on $X$ and $Y$ respectively.
For a natural number $N$ we define new metrics $d_N$ and $d^\prime_N$ on $X$ and $Y$ by 
\[  d_N(x_1, x_2) = \max_{0\leq n < N} d\left(T^n x_1, T^n x_2\right), \quad
     d^\prime_N(y_1, y_2) = \max_{0\leq n< N} d^\prime\left(S^n y_1, S^n y_2\right). \] 
For $\varepsilon>0$ and a subset $U \subset X$ we define 
\[ \#\left(U, d_N, \varepsilon\right) = \min\left\{n\geq 1\middle|
     \parbox{3in}{\centering $\exists$ open subsets $U_1, \dots, U_n$ of $X$ with $U \subset U_1\cup \dots \cup U_n$ and 
     $\diam(U_k, d_N) < \varepsilon$ for all $1\leq k\leq n$}
    \right\}. \]
When $U$ is empty, we set $\#\left(U, d_N, \varepsilon\right)=0$.

Recall $0\leq w \leq 1$. We set 
\begin{equation*} 
   \begin{split}
      & \#^w\left(\pi, X, d_N,d^\prime_N, \varepsilon\right) \\
      & = \min\left\{\sum_{k=1}^n \left(\#\left(\pi^{-1}(V_k), d_N, \varepsilon\right)\right)^w\middle|
           \parbox{3in}{\centering $Y=V_1\cup \dots \cup V_n$ is an open cover with $\diam \left(V_k, d^\prime_N\right) < \varepsilon$
          for all $1\leq k\leq n$}\right\}. 
  \end{split}    
\end{equation*}      
This quantity is sub-multiplicative in $N$ and monotone in $\varepsilon$.
We define the \textbf{$w$-weighted topological entropy} of the factor map $\pi:(X, T)\to (Y, S)$ by
\[ \htop^w(\pi, X, T) = \lim_{\varepsilon\to 0} \left(\lim_{N\to \infty} 
     \frac{\log \#^w\left(\pi, X, d_N,d^\prime_N, \varepsilon\right)}{N}\right).  \]
The value of $\htop^w(\pi,X, T)$ is independent of the choices of metrics $d$ and $d^\prime$.
So it provides a topological invariant.

The weighted topological entropy satisfies the following \textit{variational principle}
(\cite[Theorem 1.4, Corollary 1.5]{Feng--Huang} and \cite[Theorem 1.3]{Tsukamoto_weighted}):
\begin{equation} \label{eq: weighted variational principle}
   \htop^w(\pi, X, T) = \sup_{\mu\in \mathscr{M}^T(X)} \left\{w h_\mu(X,T) + (1-w)h_{\pi_*\mu}(Y,S)\right\}.
\end{equation}   
Here $\mathscr{M}^T(X)$ is the set of $T$-invariant Borel probability measures on $X$, and 
$h_\mu(X,T)$ and $h_{\pi_*\mu}(Y,S)$ are the Kolmogorov--Sinai entropy of the measure preserving systems
$(X, T, \mu)$ and $(Y, S, \pi_*\mu)$ respectively.

In the next subsection we need to use the weighted topological entropy for a factor map between symbolic dynamical systems.
In the case of symbolic dynamics, the above formulation of the weighted topological entropy is essentially the same with the one
given in \cite[Theorem 1.1]{Barral--Feng_arXiv} and \cite[Theorem 3.1]{Barral--Feng}.

Let $A$ and $B$ be finite sets, and let $\left((A\times B)^{\mathbb{N}}, \sigma\right)$ and 
$\left(B^{\mathbb{N}},\sigma\right)$ be the full-shifts on the alphabets
$A\times B$ and $B$ respectively.
Let $\pi:(A\times B)^{\mathbb{N}}\to B^{\mathbb{N}}$ be the natural projection.
Let $\Omega\subset (A\times B)^{\mathbb{N}}$ be a subshift (shift-invariant closed subset).
Set $\Omega^\prime := \pi(\Omega)\subset B^{\mathbb{N}}$.
For a natural number $N$, we define $\Omega|_N\subset (A\times B)^N$ and 
$\Omega^\prime|_N\subset B^N$ as the images of the 
projections of $\Omega$ and $\Omega^\prime$ to the first $N$ coordinates, 
respectively.
We denote by $\pi_N:\Omega|_N\to \Omega^\prime|_N$ the natural projection map.

\begin{lemma} \label{lemma: weighted topological entropy for symbolic dynamics}
In the above setting, for $0\leq w \leq 1$, the weighted topological entropy of 
the factor map $\pi:(\Omega, \sigma) \to (\Omega^\prime, \sigma)$ is given by 
\[  \htop^w\left(\pi, \Omega, \sigma\right) = \lim_{N\to \infty} \frac{1}{N}
     \log \sum_{v\in \Omega^\prime|_N} \left|\pi_N^{-1}(v)\right|^w. \]
Here $\left|\pi_N^{-1}(v)\right|$ is the cardinality of $\pi_N^{-1}(v)\subset \Omega|_N$.
\end{lemma}

Notice that $\sum_{v\in \Omega^\prime|_N} \left|\pi_N^{-1}(v)\right|^w$ is sub-multiplicative in $N$.
So the above limit exists.

\begin{proof}
It is immediate to check 
\[  \htop^w\left(\pi, \Omega, \sigma\right) \leq  \lim_{N\to \infty} \frac{1}{N}
     \log \sum_{v\in \Omega^\prime|_N} \left|\pi_N^{-1}(v)\right|^w. \]
Here we prove the reverse inequality.
We define metrics $d$ and $d^\prime$ on $(A\times B)^{\mathbb{N}}$ and $B^{\mathbb{N}}$ by 
\[ d\left((x, y), (x^\prime, y^\prime)\right) = 2^{-\min\{n\mid (x_n, y_n)\neq (x^\prime_n, y^\prime_n)\}}, \quad 
    d^\prime(y, y^\prime) = 2^{-\min\{n\mid y_n \neq y^\prime_n\}}.  \]
Let $N$ be a natural number.
For subsets $U\subset (A\times B)^{\mathbb{N}}$ and $V\subset B^{\mathbb{N}}$, we denote by 
$U|_N\subset (A\times B)^N$ and $V|_N \subset B^N$ the projections to the first $N$-coordinates.
If $\diam(U, d_N) < 1$ then $U|_N$ is a singleton or empty.
Similarly, if $\diam(V, d^\prime_N) < 1$ then so is $V|_N$.

Let $0<\varepsilon<1$.
Let $\Omega^\prime = V_1\cup \dots \cup V_n$ be an open cover with 
$\diam(V_k, d^\prime_N) < \varepsilon$ and 
\[  \#^w\left(\pi, \Omega, d_N,d^\prime_N, \varepsilon\right) 
    = \sum_{k=1}^n \left(\#\left(\pi^{-1}(V_k), d_N, \varepsilon\right)\right)^w. \]
We also assume that any $V_k$ is not empty. (Hence $V_k|_N$ is a singleton.)
Set $t_k = \#\left(\pi^{-1}(V_k), d_N, \varepsilon\right)$.

For each $V_k$ we take an open cover 
$\pi^{-1}(V_k) = U_{k1}\cup U_{k2}\cup \dots\cup U_{k t_k}$ with 
$\diam\left(U_{k l}, d_N\right) < \varepsilon$ for all $1\leq l \leq t_k$.
Any $U_{kl}$ is not empty. So $U_{kl}|_N$ is a singleton.
We have 
\begin{align*}
    \#^w\left(\pi, \Omega, d_N,d^\prime_N, \varepsilon\right) &= \sum_{k=1}^n t_k^w  \\ 
    & =  \sum_{v\in \Omega^\prime|_N}\left(\sum_{k: V_k|_N = \{v\}} t_k^w\right).  
\end{align*}    
Here, in the sum $\sum_{k: V_k|_N = \{v\}}$ in the second line, $k$ runs over all index $k\in [1, n]$ satisfying 
$V_k|_N = \{v\}$.
Since we assume $0\leq w\leq 1$,
\[  \sum_{k: V_k|_N = \{v\}} t_k^w  \geq \left(\sum_{k: V_k|_N = \{v\}} t_k\right)^w \quad 
     \left(\text{by } x^w+y^w\geq (x+y)^w \text{ for } x, y\geq 0\right). \]
For $v\in \Omega^\prime|_N$
\[ \pi_N^{-1}(v)   = \bigcup_{k: V_k|_N = \{v\}} \left(\pi^{-1}(V_k)\right)|_N 
     =  \bigcup_{k: V_k|_N = \{v\}}  \bigcup_{l=1}^{t_k} \left(U_{kl}|_N\right).  \]
Since every $U_{kl}|_N$ is a singleton, 
\[  \left|\pi_N^{-1}(v)\right| \leq \sum_{k: V_k|_N = \{v\}} t_k. \]
Therefore 
\begin{align*}
    \sum_{v\in \Omega^\prime|_N} \left|\pi_N^{-1}(v)\right|^w 
     & \leq  \sum_{v\in \Omega^\prime|_N} \left(\sum_{k: V_k|_N = \{v\}} t_k\right)^w \\
      & \leq  \sum_{v\in \Omega^\prime|_N}  \left( \sum_{k: V_k|_N = \{v\}} t_k^w\right) \\
      &=  \#^w\left(\pi, \Omega, d_N,d^\prime_N, \varepsilon\right).
\end{align*}      
Thus 
\[   \lim_{N\to \infty} \frac{1}{N}\log \sum_{v\in \Omega^\prime|_N} \left|\pi_N^{-1}(v)\right|^w \leq 
      \lim_{N\to \infty} \frac{1}{N}\log \#^w\left(\pi, \Omega, d_N,d^\prime_N, \varepsilon\right) 
       \leq \htop^w\left(\pi, \Omega, \sigma\right). \]
\end{proof}

\begin{example} \label{example: weighted topological entropy}
Let $R\subset A\times B$ be a nonempty subset. 
Set $\Omega = R^{\mathbb{N}}$. This is a subshift of $\left(A\times B\right)^{\mathbb{N}}$.
For each $v\in B$ we denote by $t(v)$ the number of $u\in A$ with $(u, v)\in R$.
Then for each $v= (v_1, \dots, v_N)\in \Omega^\prime|_N$ we have 
\[  \left|\pi_N^{-1}(v)\right| = t(v_1) \cdots t(v_N). \]
By Lemma \ref{lemma: weighted topological entropy for symbolic dynamics}
   \begin{align*}
     \htop^w\left(\pi, \Omega, \sigma\right)  
    & = \lim_{N\to \infty} \frac{1}{N} \log \sum_{v \in \Omega^\prime|_N} \left|\pi^{-1}_N(v)\right|^w \\
    & = \lim_{N\to \infty}  \frac{1}{N} \log \left(\sum_{v\in B} t(v)^w\right)^N \\
    & = \log \sum_{v\in B} t(v)^w.
   \end{align*}
\end{example}

Readers can find in \cite{Kenyon--Peres} many wonderful calculations of 
$\htop^w(\pi, \Omega, \sigma)$ in the case that $\Omega$ is a subshift of finite type or a 
sofic subshift\footnote{Kenyon--Peres \cite{Kenyon--Peres} did not use the terminologies of weighted topological 
entropy, but we can interpret their results in terms of $\htop^w\left(\pi, \Omega, \sigma\right)$.}.

\subsection{Main result for infinite dimensional carpets} \label{subsection: main result for infinite dimensional carpets}

Let $a\geq b\geq 2$ be two natural numbers and set 
\[  A = \{0,1,2,\dots, a-1\}, \quad B = \{0,1,2,\dots, b-1\}. \]
Let $(A\times B)^\mathbb{N}$ be the one-sided full-shift on the alphabet $A\times B$
with the shift map $\sigma: (A\times B)^{\mathbb{N}}\to (A\times B)^{\mathbb{N}}$.
Let $\pi: (A\times B)^{\mathbb{N}}\to B^{\mathbb{N}}$ be the natural projection.
We also denote by $\sigma:B^{\mathbb{N}}\to B^{\mathbb{N}}$ the shift map on $B^{\mathbb{N}}$.

Let $[0,1]^\mathbb{N} = [0,1] \times [0,1] \times [0,1]\times \dots$ be the infinite dimensional cube, and consider the product
$[0,1]^{\mathbb{N}}\times [0,1]^{\mathbb{N}}$ with a metric $d$ defined by 
\begin{equation}  \label{eq: metric on infinite dimensional carpets}
   d\left((x, y), (x^\prime, y^\prime)\right)  = \sum_{n=1}^\infty 2^{-n} \max\left(|x_n-x^\prime_n|, |y_n-y^\prime_n|\right), 
\end{equation}   
where $x = (x_n)_{n\in \mathbb{N}}, y= (y_n)_{n\in \mathbb{N}}$ and 
$x^\prime = (x^\prime_n)_{n\in \mathbb{N}}, y^\prime = (y^\prime_n)_{n\in \mathbb{N}}$ are points in $[0,1]^\mathbb{N}$.
We define a shift map on $[0,1]^{\mathbb{N}}\times [0,1]^{\mathbb{N}}$ 
(also denoted by $\sigma: [0,1]^{\mathbb{N}}\times [0,1]^{\mathbb{N}} \to [0,1]^{\mathbb{N}}\times [0,1]^{\mathbb{N}}$) by 
\[  \sigma\left((x_n)_{n\in \mathbb{N}}, (y_n)_{n\in \mathbb{N}}\right) 
    = \left((x_{n+1})_{n\in \mathbb{N}}, (y_{n+1})_{n\in \mathbb{N}}\right). \]

Let $\Omega\subset (A\times B)^{\mathbb{N}}$ be a subshift, namely a closed subset satisfying $\sigma(\Omega) \subset \Omega$.
We assume $\Omega \neq \emptyset$.
We define a \textbf{carpet system} $X_\Omega \subset [0,1]^{\mathbb{N}}\times [0,1]^{\mathbb{N}}$ by 
\[  X_\Omega = \left\{\left(\sum_{m=1}^\infty \frac{x_m}{a^m}, \sum_{m=1}^\infty \frac{y_m}{b^m}\right) 
      \in [0,1]^{\mathbb{N}} \times [0,1]^{\mathbb{N}} \middle|\, 
     (x_m, y_m)\in \Omega \text{ for all $m \geq 1$}\right\}. \]
Here $x_m\in A^{\mathbb{N}}\subset \ell^\infty$ and we consider the summation 
$\sum_{m=1}^\infty \frac{x_m}{a^m}$ in $\ell^\infty$. 
Then $\sum_{m=1}^\infty \frac{x_m}{a^m} \in [0,1]^{\mathbb{N}}$.
Similarly for the term $\sum_{m=1}^\infty \frac{y_m}{b^m}$.

$(X_\Omega, \sigma)$ is a subsystem of 
$\left([0,1]^{\mathbb{N}}\times [0,1]^{\mathbb{N}}, \sigma\right)$.
We are interested in its mean Hausdorff dimension and metric mean dimension with respect to the metric 
(\ref{eq: metric on infinite dimensional carpets}).
Recall that $\pi:(A\times B)^{\mathbb{N}}\to B^{\mathbb{N}}$ is the natural projection.
We consider its restriction to $\Omega$ and also denote it by $\pi:\Omega\to \pi(\Omega)$.
Set $\Omega^\prime = \pi(\Omega)$. 
This is a subshift of $B^{\mathbb{N}}$. 

\begin{theorem} \label{theorem: mean Hausdorff dimension and metric mean dimension of carpet systems}
In the above setting, the mean Hausdorff dimension and metric mean dimension of 
$(X_\Omega, \sigma, d)$ are given by 
\begin{align}
   \mdimh\left(X_\Omega, \sigma, d\right) & =  \frac{\htop^{\log_a b}\left(\pi, \Omega, \sigma\right)}{\log b}, 
    \label{eq: mean Hausdorff dimension of carpet systems} \\
   \mdimm\left(X_\Omega, \sigma, d\right) & = \frac{\htop\left(\Omega, \sigma\right)}{\log a} + 
      \left(\frac{1}{\log b} - \frac{1}{\log a}\right) \htop\left(\Omega^\prime, \sigma\right).
   \label{eq: metric mean dimension of carpet systems}
\end{align}
Here $\htop^{\log_a b}\left(\pi, \Omega, \sigma\right)$ is the weighted topological entropy of the factor map 
$\pi:(\Omega, \sigma) \to \left(\Omega^\prime, \sigma\right)$ with the weight $\log_a b = \frac{\log b}{\log a}$.
\end{theorem}

By the variational principle (\ref{eq: weighted variational principle}) for weighted topological entropy, 
the above formula of the mean Hausdorff dimension can be also expressed as
\begin{align*}
   \mdimh\left(X_\Omega, \sigma, d\right) &= \frac{1}{\log b} 
     \sup_{\mu\in \mathscr{M}^\sigma(\Omega)} \left\{(\log_a b) h_{\mu}(\Omega, \sigma) 
      + (1-\log_a b) h_{\pi_*\mu}\left(\Omega^\prime, \sigma\right)\right\}, \\
      &= \sup_{\mu\in \mathscr{M}^\sigma(\Omega)} \left\{\frac{h_\mu(\Omega, \sigma)}{\log a} + 
           \left(\frac{1}{\log b} - \frac{1}{\log a}\right) h_{\pi_*\mu}\left(\Omega^\prime,\sigma\right)\right\}.
\end{align*}     
On the other hand, by the standard variational principle for topological entropy
\[ \mdimm\left(X_\Omega, \sigma, d\right) = 
    \sup_{\mu\in \mathscr{M}^\sigma(\Omega)} \frac{h_\mu(\Omega, \sigma)}{\log a} + 
    \left(\frac{1}{\log b} - \frac{1}{\log a}\right) \sup_{\nu\in \mathscr{M}^\sigma\left(\Omega^\prime\right)}
    h_\nu\left(\Omega^\prime, \sigma\right). \]
So the difference lies in whether we take supremum simultaneously or separately for the two terms 
$h_\mu(\Omega, \sigma)$ and $h_\nu\left(\Omega^\prime, \sigma\right)$.

\begin{example}
Let $a=3$ and $b=2$. Then $A= \{0,1,2\}$ and $B = \{0,1\}$.
Set 
\[  R = \{(0,0), (1,1), (2,0)\}  \subset A\times B. \]
Define $\Omega = R^{\mathbb{N}}$. 
Then $\Omega^\prime = B^{\mathbb{N}}$.
By Example \ref{example: weighted topological entropy}
\[   \htop^{\log_3 2}(\pi, \Omega, \sigma) = \log \left(1 + 2^{\log_3 2}\right). \]
It is also easy to see 
\[  \htop(\Omega, \sigma) = \log 3, \quad 
     \htop\left(\Omega^\prime, \sigma\right) = \log 2. \]
Therefore by 
Theorem \ref{theorem: mean Hausdorff dimension and metric mean dimension of carpet systems}
\begin{align*}
   \mdimh\left(X_\Omega, \sigma, d\right) & = \log_2\left(1+2^{\log_3 2}\right) = 1.3496838201\dots,  \\
   \mdimm\left(X_\Omega, \sigma, d\right) & = 1 + 
    \left(\frac{1}{\log 2} - \frac{1}{\log 3}\right) \log 2 \\
   & = 2- \log_3 2 = 1.3690702464\dots.
\end{align*}
As we already mentioned at the end of \S \ref{subsection: weighted topological entropy}, readers
can find many interesting calculations of $\htop^{\log_a b}\left(\pi, \Omega, \sigma\right)$ for sofic subshifts $\Omega$
in \cite{Kenyon--Peres}.
\end{example}

\subsection{Calculation of metric mean dimension} \label{subsection: calculation of metric mean dimension of carpets}

Here we prove the formula (\ref{eq: metric mean dimension of carpet systems}) 
of the metric mean dimension of the carpet system $(X_\Omega, \sigma, d)$.
We use the notations of Theorem 
\ref{theorem: mean Hausdorff dimension and metric mean dimension of carpet systems}.
Set 
\[   w = \log_a b = \frac{\log b}{\log a}. \]
Since $a\geq b \geq 2$, we have $0<w\leq 1$.
For a natural number $N$ we denote by $\Omega|_N$ and $\Omega^\prime|_N$ the images of $\Omega$ and $\Omega^\prime$ 
under the projections
\begin{align*}
  (A\times B)^{\mathbb{N}} & \to  A^N \times B^N,  \quad 
     \left((u_n)_{n\in \mathbb{N}}, (v_n)_{n\in \mathbb{N}}\right) \mapsto \left((u_1, \dots, u_N), (v_1, \dots, v_N)\right), \\
  B^{\mathbb{N}} & \to B^N,  \quad    (v_n)_{n\in \mathbb{N}} \mapsto (v_1, \dots, v_N).
\end{align*}
We also define $X_\Omega |_N$ as the image of $X_\Omega$ under the projection 
\[  [0,1]^{\mathbb{N}}\times [0,1]^{\mathbb{N}} \to [0,1]^N \times [0,1]^N, \quad 
     \left((x_n)_{n\in \mathbb{N}}, (y_n)_{n\in \mathbb{N}}\right) \mapsto \left((x_1, \dots, x_N), (y_1, \dots, y_N)\right). \]
We have 
\[  X_\Omega|_N = \left\{\left(\sum_{m=1}^\infty \frac{x_m}{a^m}, \sum_{m=1}^\infty \frac{y_m}{b^m}\right)
     \in [0,1]^N\times [0,1]^N\middle|\, (x_m,y_m)\in \Omega|_N \text{ for all } m\geq 1\right\}. \]

Let $(x, y) \in \left(\Omega|_N\right)^{\mathbb{N}}$ where $x = (x_m)_{m\in \mathbb{N}}$ and $y = (y_m)_{m\in \mathbb{N}}$ with 
$x_m\in A^N$, $y_m\in B^N$ and $(x_m, y_m) \in \Omega|_N$.
For natural numbers $N$ and $M$, we define a subset $Q_{N,M}(x,y) \subset X_\Omega |_N$ by
\begin{equation}  \label{eq: definition of Q_NM}
    Q_{N, M}(x, y)  =
   \left\{\left(\sum_{m=1}^\infty \frac{x^\prime_m}{a^m}, \sum_{m=1}^\infty \frac{y^\prime_m}{b^m}\right)  \middle|\, 
    \parbox{2.5in}{\centering $(x^\prime_m, y^\prime_m) \in \Omega|_N$ for all $m\geq 1$ with \\
    $x^\prime_m = x_m \> (1\leq m \leq \lfloor wM \rfloor)$ and \\ $y^\prime_m = y_m \> (1\leq m \leq M)$} \right\}. 
\end{equation}   
Here $\lfloor wM\rfloor$ is the largest integer not greater than $wM$. 
The set $Q_{N,M}(x,y)$ depends only on the coordinates $x_1, \dots, x_{\lfloor wM\rfloor}, y_1, \dots, y_M$. So we also denote it by 
\[  Q_{N,M}(x_1,\dots, x_{\lfloor wM \rfloor}, y_1, \dots, y_M) \quad \left( =  Q_{N, M}(x, y)\right). \]

From $w = \log_a b$, we have $a^{-wM} = b^{-M}$ and
\[  b^{-M} \leq a^{-\lfloor wM\rfloor} < a b^{-M}. \]
We have
\begin{equation}   \label{eq: diameter of Q_NM}
    \diam\left(Q_{N,M}(x,y), \norm{\cdot}_\infty\right) \leq \max\left(a^{-\lfloor wM\rfloor}, b^{-M}\right)  = a^{-\lfloor wM\rfloor} 
     < a b^{-M}. 
\end{equation}     
Here $\norm{\cdot}_\infty$ is the $\ell^\infty$-distance (i.e. the distance defined by the $\ell^\infty$-norm) 
on $X_\Omega|_N \subset \mathbb{R}^{2N}$.

\begin{lemma} \label{lemma: projection of X_Omega to the first N coordinates}
Let $\varepsilon$ be a positive number. For any $N\geq 1$ we have 
\[  \#\left(X_\Omega|_N, \norm{\cdot}_\infty, \varepsilon\right) \leq \#\left(X_\Omega, d_N, \varepsilon\right). \]
Let $L$ be a natural number satisfying $\sum_{n>L} 2^{-n} < \varepsilon/2$. Then 
\[   \#\left(X_\Omega, d_N, \varepsilon\right) \leq \#\left(X_\Omega|_{N+L}, \norm{\cdot}_\infty, \frac{\varepsilon}{2}\right). \]
\end{lemma}

\begin{proof}
The first inequality follows from 
\[  \norm{x|_N -y|_N}_\infty \leq d_N(x,y) \quad (x,y\in X_\Omega), \]
where $x|_N$ denotes the projection of $x$ to $X_\Omega|_N$.
The second inequality follows from 
\[   d_N(x,y) < \norm{x|_{N+L} - y|_{N+L}}_\infty + \frac{\varepsilon}{2}. \]
\end{proof}

\begin{lemma}  \label{lemma: covering number of X_Omega|_N}
For any natural numbers $N$ and $M$ we have 
   \begin{align*}
      \#\left(X_\Omega|_N, \norm{\cdot}_\infty, a b^{-M}\right) 
      & \leq \left|\Omega|_N\right|^{\lfloor wM\rfloor}\cdot \left|\Omega^\prime|_N\right|^{M-\lfloor wM\rfloor}, \\
      \#\left(X_\Omega|_N, \norm{\cdot}_\infty, b^{-M}\right) & \geq 
      \left|\Omega|_N\right|^{\lfloor wM\rfloor}\cdot \left|\Omega^\prime|_N\right|^{M-\lfloor wM\rfloor}.
   \end{align*}
Here $\left|\Omega|_N\right|$ and $\left|\Omega^\prime|_N\right|$ 
denote the cardinalities of $\Omega|_N$ and $\Omega^\prime|_N$
respectively.
\end{lemma}

\begin{proof}
Recall that $\diam\left(Q_{N,M}(x,y), \norm{\cdot}_\infty\right) < a b^{-M}$.
The first inequality follows from 
\begin{equation*}
 X_\Omega|_N = \bigcup\left\{Q_{N,M}(x_1,\dots,x_{\lfloor wM\rfloor}, y_1, \dots, y_M)\middle|\, 
  \parbox{2.6in}{\centering $(x_m,y_m)\in \Omega|_N$ for $1\leq m \leq \lfloor wM\rfloor$ \\
                    $y_m\in \Omega^\prime|_N$ for $\lfloor wM\rfloor +1 \leq m \leq M$}\right\}.
\end{equation*}

Next we consider the second inequality.
We fix a point $(\xi, \eta)\in \Omega|_N$.
For each $v \in \Omega^\prime|_N$ we pick up $s(v)\in A^N$ satisfying $\left(s(v),v\right) \in \Omega|_N$.
For $(x_m, y_m) \in \Omega|_N$ $(1\leq m \leq \lfloor wM\rfloor)$ and 
$y_m\in \Omega^\prime|_N$ $(\lfloor wM\rfloor +1 \leq m \leq M)$, we set 
\begin{align*}
 &  p(x_1, \dots, x_{\lfloor wM\rfloor}, y_1, \dots, y_M) \\
   & =   \left(\sum_{m=1}^{\lfloor wM\rfloor} \frac{x_m}{a^m} + \sum_{m=\lfloor wM\rfloor +1}^M \frac{s(y_m)}{a^m} 
           + \sum_{m=M+1}^\infty \frac{\xi}{a^m}, \sum_{m=1}^M \frac{y_m}{b^m} + \sum_{m=M+1}^\infty \frac{\eta}{b^m}\right). 
\end{align*}           
This is a point in $X_\Omega|_N$.
For $(x_1, \dots, x_{\lfloor wM\rfloor}, y_1, \dots, y_M) \neq
(x^\prime_1, \dots, x^\prime_{\lfloor wM\rfloor}, y^\prime_1, \dots, y^\prime_M)$, we have 
\begin{align*}
    \norm{p(x_1, \dots, x_{\lfloor wM\rfloor}, y_1, \dots, y_M) 
    - p(x^\prime_1, \dots, x^\prime_{\lfloor wM\rfloor}, y^\prime_1, \dots, y^\prime_M)}_\infty
    & \geq \min\left(a^{-\lfloor wM\rfloor}, b^{-M}\right) \\
    & = b^{-M}. 
\end{align*}    
Hence the set 
\[ \left\{ p(x_1, \dots, x_{\lfloor wM\rfloor}, y_1, \dots, y_M)\middle|\,  
     \parbox{2.6in}{\centering $(x_m,y_m)\in \Omega|_N$ for $1\leq m \leq \lfloor wM\rfloor$ \\
                    $y_m\in \Omega^\prime|_N$ for $\lfloor wM\rfloor +1 \leq m \leq M$}\right\} \]
is $b^{-M}$-separated with respect to the $\ell^\infty$-distance.
The second inequality follows from this.
\end{proof}

From Lemmas \ref{lemma: projection of X_Omega to the first N coordinates} and 
\ref{lemma: covering number of X_Omega|_N} we can calculate the metric mean dimension of 
$(X_\Omega, \sigma, d)$:
\begin{align*}
   \mdimm\left(X_\Omega, \sigma, d\right) & = 
   \lim_{M\to \infty} \left\{\lim_{N\to \infty} 
   \frac{\log \left(\left|\Omega|_N\right|^{\lfloor wM\rfloor}\cdot 
   \left|\Omega^\prime|_N\right|^{M-\lfloor wM\rfloor}\right)}{NM\log b}\right\} \\
    &= \lim_{M\to \infty} \left\{\frac{1}{M\log b} \lim_{N\to \infty}
    \left(\frac{\lfloor wM\rfloor \log \left|\Omega|_N\right| 
    + \left(M-\lfloor wM\rfloor\right)\log \left|\Omega^\prime|_N\right|}{N}\right)\right\}  \\
    & = \lim_{M\to \infty} \frac{\lfloor wM\rfloor \htop(\Omega, \sigma) 
     + \left(M-\lfloor wM\rfloor\right) \htop\left(\Omega^\prime, \sigma\right)}{M\log b} \\
   & = \frac{w}{\log b} \htop(\Omega, \sigma) + \frac{1-w}{\log b} \htop\left(\Omega^\prime, \sigma\right)  \\
   & =  \frac{\htop\left(\Omega, \sigma\right)}{\log a} + 
      \left(\frac{1}{\log b} - \frac{1}{\log a}\right) \htop\left(\Omega^\prime, \sigma\right), \quad 
      \left(w = \log_a b = \frac{\log b}{\log a}\right).
\end{align*}    
This proves the formula (\ref{eq: metric mean dimension of carpet systems}) in 
Theorem \ref{theorem: mean Hausdorff dimension and metric mean dimension of carpet systems}.

\subsection{Preparations}  \label{subsection: preparations for infinite dimensional carpets}

The calculation in the last subsection
was rather straightforward.
It is more difficult to calculate the mean Hausdorff dimension.
This subsection is a preparation for it.

First we prepare a general result on metric geometry.
This result will be also used in the Appendix.

\begin{lemma} \label{lemma: upper bound on epsilon Hausdorff dimension}
Let $c, \varepsilon, s$ be positive numbers.
Let $(X, d)$ be a compact metric space with a Borel probability measure $\mu$.
Suppose:
\begin{itemize}
   \item $6\varepsilon^c <1$.
   \item For any $x\in X$ there exists a Borel subset $A\subset X$ satisfying $x\in A$ and
   \[   0< \diam\, A < \frac{\varepsilon}{6}, \quad \mu(A) \geq \left(\diam\, A\right)^s. \]
\end{itemize}
Then $\dimh(X, d, \varepsilon) \leq (1+c) s$.
\end{lemma}

\begin{proof}
For any $x\in X$ there exists $A_x\subset X$ with $x\in A_x$ satisfying
$0<\diam\, A_x < \varepsilon/6$ and $\mu(A_x) \geq \left(\diam\, A_x\right)^s$. 
The set $A_x$ is contained in $B\left(x, \diam\, A_x\right)$ (the closed ball of radius $\diam\, A_x$ centered at $x$).
By the compactness, there are $x_1, \dots, x_n\in X$ satisfying (let $r_k := \diam\, A_{x_k}$)
\[   X = \bigcup_{k=1}^n B(x_k, r_k), \quad \mu\left(B(x_k, r_k)\right) \geq \mu(A_{x_k}) \geq r_k^s. \]
By the finite Vitali covering lemma \cite[Lemma 2.27]{Einsiedler--Ward}, 
we can pick $x_{k_1}, \dots, x_{k_m}$ such that 
\begin{itemize}
    \item $B(x_{k_1}, r_{k_1}), \dots, B(x_{k_m}, r_{k_m})$ are mutually disjoint.
    \item $X= \bigcup_{j=1}^m B(x_{k_j}, 3 r_{k_j})$.
\end{itemize}
Then $\diam B(x_{k_j}, 3 r_{k_j}) \leq 6 r_{k_j} < \varepsilon$ and 
\begin{align*}
    \sum_{j=1}^m \left(6r_{k_j}\right)^{(1+c)s} & = \sum_{j=1}^m \left(6 r_{k_j}\right)^{cs} \cdot \left(6r_{k_j}\right)^{s} \\
      & < \sum_{j=1}^m \varepsilon^{cs} \cdot \left(6r_{k_j}\right)^{s} \\
      &= \sum_{j=1}^m \left(6\varepsilon^c\right)^s \cdot r_{k_j}^s \\
      & < \sum_{j=1}^m r_{k_j}^s \quad (\text{by $6\varepsilon^c <1$}) \\
      & \leq \sum_{j=1}^m \mu\left(B(x_{k_j}, r_{k_j})\right) \quad 
       (\text{by $\mu\left(B(x_k, r_k)\right)  \geq r_k^s$}) \\
      & \leq 1 \quad (\text{since $B(x_{k_1}, r_{k_1}), \dots, B(x_{k_m}, r_{k_m})$ are mutually disjoint}).
\end{align*}      
This shows $\dimh(X, d, \varepsilon) \leq (1+c)s$.
\end{proof}

We also need some basic results of probability theory.
For a random variable $\xi$ we denote its mean and variance
by $\mathbb{E}(\xi)$ and 
$\mathbb{V}(\xi) := \mathbb{E}[(\xi-\mathbb{E}\xi)^2]$ respectively.
Recall the Kolmogorov maximal inequality 
(see, e.g. \cite[Theorem 2.5.2]{Durrett}, \cite[Chapter IX, \S 7]{Feller}):

\begin{theorem}[Kolmogorov maximal inequality]
Let $\xi_1, \dots, \xi_n$ be independent random variables with $\mathbb{E}(\xi_k) = 0$ and $\mathbb{V}(\xi_k) < \infty$ 
for all $1\leq k \leq n$.  Set $S_k = \xi_1+\dots+\xi_k$. Then for any positive number $h$
\[   \mathbb{P}\left(\max_{1\leq k \leq n} |S_k| \geq h \right) \leq \frac{\mathbb{V}(S_n)}{h^2}. \]
\end{theorem}

The next lemma shows simple consequences of this theorem.

\begin{lemma} \label{lemma: law of large numbers}
   \begin{enumerate}
    \item  Let $\xi_1, \xi_2, \xi_3, \dots$ be independent and identically distributed (i.i.d.) 
              random variables with $\xi_1\in L^2$. 
              Then for any positive number $\delta$ and natural number $m$
    \[  \mathbb{P}\left(\sup_{n \geq m} \left|\frac{\xi_1+\dots+\xi_n}{n} - \mathbb{E}\xi_1\right|>\delta\right) 
         \leq \frac{4\mathbb{V}(\xi_1)}{\delta^2 m}. \]
    \item  For any positive numbers $C$ and $\delta$ there exists a natural number $m_0 = m_0(C,\delta)$ for which the following 
    statement holds true: If $\xi_1, \xi_2, \xi_3, \dots$ be i.i.d. random variables with $\xi_1\in L^2$ and $\mathbb{V}\xi_1 \leq C$,
    then
    \[   \mathbb{P}\left(\sup_{n\geq m_0}  \left|\frac{\xi_1+\dots+\xi_n}{n} - \mathbb{E}\xi_1\right| \leq \delta\right)  \geq \frac{3}{4}. \]
    \item For any $0<w\leq 1$ and any positive numbers $C$ and $\delta$, there exists a natural number 
    $m_1 = m_1(w, C, \delta)$ for which the following statement holds true:
    If  $\xi_1, \xi_2, \xi_3, \dots$ be i.i.d. random variables with $\xi_1\in L^2$, $\left|\mathbb{E}\xi_1\right|\leq C$ 
    and $\mathbb{V}\xi_1 \leq C$ then
    \[  \mathbb{P}\left(\sup_{n\geq m_1} \left|\frac{\xi_1+\dots+\xi_n}{n} - \frac{\xi_1+\dots+\xi_{\lfloor wn \rfloor}}{wn}\right| \leq \delta\right)
         \geq \frac{1}{2}. \]
   \end{enumerate}
\end{lemma}

\begin{proof}
(1) Let $\sigma^2 = \mathbb{V}(\xi_1)$.
Set $\eta_n = \xi_n-\mathbb{E}\xi_n$ and $S_n = \eta_1+\dots+\eta_n$.
We have $\mathbb{V}(S_n) = n \sigma^2$.
By the Kolmogorov maximal inequality, for any $h>0$ and natural number $n$
\[  \mathbb{P}\left(\max_{n\leq k \leq 2n} |S_k| \geq h\right)
     \leq \mathbb{P}\left(\max_{1\leq k \leq 2n} |S_k|\geq h\right) \leq \frac{2n \sigma^2}{h^2}. \]
Letting $h = \delta n$, we get 
\[  \mathbb{P}\left(\max_{n\leq k \leq 2n} |S_k|\geq \delta n\right) \leq \frac{2\sigma^2}{\delta^2 n}. \]
We have 
\[  \left\{\max_{n \leq k \leq 2n} \left|\frac{S_k}{k}\right|\geq \delta\right\} \subset 
     \left\{\max_{n\leq k \leq 2n} |S_k|\geq \delta n\right\}. \]
Hence 
\[  \mathbb{P}\left(\max_{n \leq k \leq 2n} \left|\frac{S_k}{k}\right|\geq \delta\right) \leq \frac{2\sigma^2}{\delta^2 n}. \]
We use this inequality for $n = m, 2m , 2^2m, 2^3m, \dots$. Then for each $i=0,1,2,\dots$
\[  \mathbb{P}\left(\max_{2^i m \leq k \leq 2^{i+1}m} \left|\frac{S_k}{k}\right|\geq \delta\right) \leq \frac{2\sigma^2}{2^i \delta^2 m}. \]
We have 
\[  \left\{\sup_{k\geq m}\left|\frac{S_k}{k}\right| > \delta\right\}
     = \bigcup_{i=0}^\infty \left\{\max_{2^i m \leq k \leq 2^{i+1}m}\left|\frac{S_k}{k}\right| > \delta\right\}. \]
Therefore 
\[  \mathbb{P}\left(\sup_{k\geq m}\left|\frac{S_k}{k}\right| > \delta\right) \leq 
     \sum_{i=0}^\infty \mathbb{P}\left(\max_{2^i m \leq k \leq 2^{i+1}m}\left|\frac{S_k}{k}\right| \geq \delta\right)
      \leq \frac{4\sigma^2}{\delta^2 m}. \]
      
(2)   Take a natural number $m_0$ with $\frac{4C}{\delta^2 m_0} \leq \frac{1}{4}$.
      By (1)
\[  \mathbb{P}\left(\sup_{n\geq m_0}  \left|\frac{\xi_1+\dots+\xi_n}{n} - \mathbb{E}\xi_1\right| > \delta\right) \leq \frac{1}{4}. \] 
Hence 
\[  \mathbb{P}\left(\sup_{n\geq m_0}  \left|\frac{\xi_1+\dots+\xi_n}{n} - \mathbb{E}\xi_1\right| \leq \delta\right)  \geq \frac{3}{4}. \]

(3) Let $m_0= m_0(C,\delta/3)$ be the natural number given in (2).
We take a natural number $m_1$ satisfying 
\[   m_1 \geq \frac{m_0}{w}, \quad m_1 \geq \frac{3C}{w\delta}. \]
By (2)
\[  \mathbb{P}\left(\sup_{n\geq m_0} \left|\frac{\xi_1+\dots+\xi_n}{n} - \mathbb{E}\xi_1\right| \leq \frac{\delta}{3} \right) \geq \frac{3}{4}. \]
For $n\geq m_1$ we have $\lfloor nw \rfloor \geq m_0$. Hence 
\[  \mathbb{P}\left(\sup_{n\geq m_1} \left|\frac{\xi_1+\dots+\xi_{\lfloor wn\rfloor}}{\lfloor wn\rfloor} 
    - \mathbb{E}\xi_1\right| \leq \frac{\delta}{3} \right) \geq \frac{3}{4}. \]
We have 
\begin{align*}
   \left|\frac{\xi_1+\dots+\xi_{\lfloor wn\rfloor}}{wn} - \mathbb{E}\xi_1\right|
     & = \left|\frac{\lfloor wn\rfloor}{wn}\left(\frac{\xi_1+\dots+\xi_{\lfloor wn\rfloor}}{\lfloor wn\rfloor} - \mathbb{E}\xi_1\right)
           + \left(\frac{\lfloor wn\rfloor}{wn}-1\right) \mathbb{E}\xi_1\right|  \\
   &\leq   \left|\frac{\xi_1+\dots+\xi_{\lfloor wn\rfloor}}{\lfloor wn\rfloor} - \mathbb{E}\xi_1\right| 
     + \frac{wn-\lfloor wn\rfloor}{wn} \left|\mathbb{E}\xi_1\right| \\
   &\leq   \left|\frac{\xi_1+\dots+\xi_{\lfloor wn\rfloor}}{\lfloor wn\rfloor} - \mathbb{E}\xi_1\right| 
              + \frac{C}{wn}.
\end{align*}
For $n\geq m_1$ we have $\frac{C}{wn} \leq \frac{C}{w m_1} \leq \frac{\delta}{3}$. Hence 
\[ \sup_{n\geq m_1} \left|\frac{\xi_1+\dots+\xi_{\lfloor wn\rfloor}}{wn} - \mathbb{E}\xi_1\right|
    \leq \sup_{n\geq m_1} \left|\frac{\xi_1+\dots+\xi_{\lfloor wn\rfloor}}{\lfloor wn\rfloor} - \mathbb{E}\xi_1\right|  + \frac{\delta}{3}. \]
Therefore 
\[ \mathbb{P}\left(\sup_{n\geq m_1} \left|\frac{\xi_1+\dots+\xi_{\lfloor wn\rfloor}}{wn} 
    - \mathbb{E}\xi_1\right| \leq \frac{2\delta}{3} \right) \geq \frac{3}{4}. \]
We have 
\begin{align*}
  & \left\{\sup_{n\geq m_1} \left|\frac{\xi_1+\dots+\xi_n}{n} - \frac{\xi_1+\dots+\xi_{\lfloor wn \rfloor}}{wn}\right| \leq \delta\right\}\\
  &   \supset \left\{\sup_{n\geq m_1} \left|\frac{\xi_1+\dots+\xi_n}{n} - \mathbb{E}\xi_1\right| \leq \frac{\delta}{3}\right\} 
     \cap \left\{\sup_{n\geq m_1} \left|\frac{\xi_1+\dots+\xi_{\lfloor wn\rfloor}}{wn} 
    - \mathbb{E}\xi_1\right| \leq \frac{2\delta}{3} \right\}.
\end{align*}     
Therefore 
\[  \mathbb{P}\left(\sup_{n\geq m_1} \left|\frac{\xi_1+\dots+\xi_n}{n} - \frac{\xi_1+\dots+\xi_{\lfloor wn \rfloor}}{wn}\right| \leq \delta\right)
     \geq \frac{1}{2}. \]
\end{proof}

We also need the following elementary lemma.

\begin{lemma} \label{lemma: elementary lemma for carpets}
Let $0<w\leq 1$ and $C\geq 0$.
Let $\{a_n\}_{n=0}^\infty$ be a sequence of real numbers with $0\leq a_n \leq C$.
For any $\delta>0$ and any natural number $m$ there exists $n$ with
\[   m \leq n \leq w^{-\lceil C/\delta\rceil}\left(m+\left\lceil \frac{C}{\delta} \right\rceil \right) + 1 \]
satisfying  
\[  a_n - a_{\lfloor wn\rfloor} \geq -\delta. \]
\end{lemma}

The number $w^{-\lceil C/\delta\rceil}\left(m+\left\lceil \frac{C}{\delta} \right\rceil \right) + 1$ looks complicated.
We will not need its precise form in the sequel.
The point is that it depends only on $w, C, \delta, m$.

\begin{proof}
Set $k = \lceil C/\delta\rceil$.
Suppose that for all integers $n\in \left[m, w^{-k}(m+k)+1\right]$ we have 
\[  a_n - a_{\lfloor wn \rfloor} < -\delta. \]
We define a non-increasing sequence $n_0 \geq  n_1 \geq  \dots \geq  n_k$ by 
\[ n_0 = \lfloor w^{-k}(k+m)+1\rfloor, \quad n_{j+1} = \lfloor w n_j\rfloor. \]
We have $n_{j+1} \geq w n_j-1$ and hence 
\[  n_{j+1}+\frac{1}{1-w} \geq w \left(n_j+\frac{1}{1-w}\right). \]
Then 
\begin{align*}
    n_k &\geq w^k \left(n_0 + \frac{1}{1-w}\right) - \frac{1}{1-w} \\
          & = w^k n_0 - \frac{1-w^k}{1-w} \\
         &= w^k n_0 - \left(1+w+w^2+\dots+w^{k-1}\right) \\
         & \geq  w^k n_0 -k \\
         & \geq w^k\left(w^{-k}(m+k)\right) -k  \quad 
\left(\text{by } n_0 = \lfloor w^{-k}(k+m)+1\rfloor \right)\\
         & =m.
\end{align*}
So all $n_0, n_1, \dots, n_k$ belong to $\left[m, w^{-k}(m+k)+1\right]$.
Therefore 
\[   a_{n_j} - a_{n_{j+1}} < -\delta \quad   \text{for all $0\leq j \leq k-1$}. \]
Since we assumed $a_{n_k} \leq C$,
\begin{align*}
    a_{n_0} &< a_{n_k} -\delta k \\
    &\leq  C -C   \quad (\text{by } k = \lceil C/\delta\rceil) \\
    & = 0.   
\end{align*}
Namely $a_{n_0} <0$. This contradicts with the assumption $a_n\geq 0$.
\end{proof}

\subsection{Calculation of mean Hausdorff dimension}  
\label{subsection: calculation of mean Hausdorff dimension of infinite dimensional carpets}

We calculate the mean Hausdorff dimension of the carpet system $(X_\Omega, \sigma, d)$ in this subsection.
We use the notations introduced in \S \ref{subsection: main result for infinite dimensional carpets} and 
\S \ref{subsection: calculation of metric mean dimension of carpets}.
Our argument is motivated by \cite{McMullen} and \cite[Section 4.2]{Bishop--Peres}.

Recall $w = \log_a b$.
Let $N$ be a natural number.
For each $v \in \Omega^\prime|_N$ we define $t_N(v)$ as the number of $u \in A^N$ with $(u,v)\in \Omega|_N$.
Set 
\[   Z_N = \sum_{v\in \Omega^\prime|_N} t_N(v)^w. \]
By Lemma \ref{lemma: weighted topological entropy for symbolic dynamics} we have 
\begin{equation}  \label{eq: formula of weighted topological entropy}
  \htop^w(\pi, \Omega, \sigma) = \lim_{N\to \infty} \frac{\log Z_N}{N}. 
\end{equation} 

For $(u,v)\in \Omega|_N$ we set 
\[  f_N(u,v) = \frac{1}{Z_N} t_N(v)^{w-1}. \]
We have 
\[  \sum_{(u,v)\in \Omega|_N} f_N(u,v) = 1. \]
So we can consider $f_N(u,v)$ as a probability measure on $\Omega|_N$.
Let $\mu_N = (f_N)^{\otimes \mathbb{N}}$ be the product of infinite copies of the measure $f_N$.
This is a probability measure defined on $\left(\Omega|_N\right)^{\mathbb{N}}$.

Let $N$ and $M$ be natural numbers.
Let $(x,y)\in \left(\Omega|_N\right)^{\mathbb{N}}$ where $x = (x_m)_{m\in \mathbb{N}}$ and 
$y = (y_m)_{m\in \mathbb{N}}$ with $x_m\in A^N$, $y_m\in B^N$ and $(x_m, y_m) \in \Omega|_N$. 
We set 
\[  P_{N,M}(x,y) = \left\{(x^\prime,y^\prime)\in \left(\Omega|_N\right)^{\mathbb{N}} \middle|\, 
    \parbox{2in}{\centering $x^\prime_m = x_m$ $(1\leq m \leq \lfloor wM\rfloor)$ \\
                     $y^\prime_m = y_m$ $(1\leq m \leq M)$}\right\}. \]                 
The set $Q_{N,M}(x,y)$ introduced in (\ref{eq: definition of Q_NM}) is the image of $P_{N,M}(x,y)$ under the map
\[  \left(\Omega|_N\right)^{\mathbb{N}} \to X_\Omega|_N, \quad 
     (x^\prime, y^\prime) \mapsto \left(\sum_{m=1}^\infty \frac{x^\prime_m}{a^m}, \sum_{m=1}^\infty \frac{y^\prime_m}{b^m}\right). \]
The key argument below is a calculation of the measure $\mu_N\left(P_{N,M}(x,y)\right)$.
Since the value of $f_N(u,v) = \frac{1}{Z_N} t_N(v)^{w-1}$ depends only on $v$, we have 
\[  \mu_N\left(P_{N,M}(x,y)\right) = \prod_{m=1}^M f_N(x_m,y_m)\cdot \prod_{m=\lfloor wM\rfloor +1}^M t_N(y_m). \]
Taking a logarithm,
\begin{align*}
   \log \mu_N\left(P_{N,M}(x,y)\right)  & = \sum_{m=1}^M \log f_N(x_m,y_m) + \sum_{m=\lfloor wM\rfloor +1}^M \log t_N(y_m) \\
  & = -M\log Z_N + (w-1)\sum_{m=1}^M \log t_N(y_m) + \sum_{m=\lfloor wM\rfloor +1}^M \log t_N(y_m) \\
  & = -M \log Z_N + w \sum_{m=1}^M \log t_N(y_m) - \sum_{m=1}^{\lfloor wM\rfloor} \log t_N(y_m)
\end{align*}
Set 
\[  S_{N,M}(x,y) = \sum_{m=1}^M \frac{\log t_N(y_m)}{N}.  \]
We also set $S_{N,0}(x,y) = 0$.
We have 
\begin{equation}  \label{eq: mu_N measure of P_NM}
   \frac{1}{NM} \log \mu_N\left(P_{N,M}(x,y)\right) = -\frac{\log Z_N}{N} 
    + w\left(\frac{S_{N,M}(x,y)}{M} - \frac{S_{N,\lfloor wM\rfloor}(x,y)}{wM}\right).
\end{equation}

\begin{lemma}  \label{lemma: lower bound on measure of P_NM}
For any positive number $\delta$ and any natural number $K$ there exists a natural number $L = L(\delta, K) \geq K$ 
for which the following 
statement holds true:
For any $N\geq 1$ and any $(x,y)\in \left(\Omega|_N\right)^{\mathbb{N}}$ there exists a natural number $M \in [K, L]$ satisfying 
\[  \frac{1}{NM} \log \mu_N\left(P_{N,M}(x,y)\right) \geq -\frac{\log Z_N}{N} - \delta. \]
\end{lemma}

\begin{proof}
From (\ref{eq: mu_N measure of P_NM})
\[   \frac{1}{NM}\log \mu_N\left(P_{N,M}(x,y)\right) \geq  -\frac{\log Z_N}{N} 
      + w\left(\frac{S_{N,M}(x,y)}{M} - \frac{S_{N,\lfloor wM\rfloor}(x,y)}{\lfloor wM\rfloor}\right). \]
We apply Lemma \ref{lemma: elementary lemma for carpets} to the sequence 
\[  a_M := \frac{S_{N,M}(x,y)}{M}. \]
We have $0\leq a_M \leq \log a$. Let 
\[  L = \left\lfloor 
     w^{-\left\lfloor \frac{\log a}{\delta}\right\rfloor} \left(K + \left\lfloor \frac{\log a}{\delta}\right\rfloor\right) + 1\right\rfloor.   \]
Then, by Lemma \ref{lemma: elementary lemma for carpets}, there exists $M\in [K, L]$ satisfying 
\[  a_M - a_{\lfloor wM\rfloor} \geq -\delta. \]
\end{proof}

\begin{lemma}  \label{lemma: upper bound on measure of P_NM}
For any $\delta>0$ there exists a natural number $m_2 = m_2(\delta)>0$ such that for any natural number $N$ 
there exists a Borel subset 
$R(\delta, N) \subset \left(\Omega|_N\right)^{\mathbb{N}}$ satisfying the following two conditions.
  \begin{itemize}
     \item $\mu_N\left(R(\delta,N)\right) \geq \frac{1}{2}$,
     \item For any $M\geq m_2$ and $(x,y)\in R(\delta,N)$ we have 
       \[  \left|\frac{1}{NM}\log \mu_N\left(P_{N,M}(x,y)\right)  -\frac{\log Z_N}{N}\right| \leq \delta. \]
  \end{itemize}
\end{lemma}

The point of the statement is that $m_2(\delta)$ is independent of $N$.

\begin{proof}
We define $\xi_m: \left(\Omega|_N\right)^{\mathbb{N}}\to \mathbb{R}$ for $m\geq 1$ by
\[  \xi_m(x, y) = \frac{\log t_N(y_m)}{N}. \] 
Notice that this depends not only on $m$ but also $N$. 
But we suppress the dependence on $N$ in our notation for simplicity.
The random variables $\xi_1, \xi_2, \xi_3, \dots$ are independent and identically distributed with respect to 
the measure $\mu_N = (f_N)^{\otimes \mathbb{N}}$.
We have $0\leq \xi_m \leq \log a$. Hence its mean and variance (with respect to $\mu_N$) are bounded by
\[  0\leq \mathbb{E}(\xi_m) \leq  \log a, \quad  \mathbb{V}(\xi_m) \leq (\log a)^2. \]
We apply Lemma \ref{lemma: law of large numbers} (3) to $\{\xi_m\}$.
Then we can find $m_2(\delta)>0$ such that the set
\[ R(\delta, N) :=\left\{(x,y)\in \left(\Omega|_N\right)^{\mathbb{N}}\middle|\,   
\sup_{M\geq m_2(\delta)} \left|\frac{\xi_1+\dots+\xi_M}{M} - \frac{\xi_1+\dots+\xi_{\lfloor wM\rfloor}}{wM}\right| \leq \delta\right\} \]
satisfies $\mu_N\left(R(\delta, N)\right) \geq \frac{1}{2}$.

We have 
\[   S_{N,M}(x,y) = \sum_{m=1}^M \xi_m. \]
From (\ref{eq: mu_N measure of P_NM})
\[  \frac{1}{NM} \log \mu_N\left(P_{N,M}(x,y)\right) + \frac{\log Z_N}{N} 
    = w\left(\frac{\xi_1+\dots+\xi_M}{M} - \frac{\xi_1+\dots+\xi_{\lfloor wM\rfloor}}{wM}\right). \]
Hence for any $(x,y)\in R(\delta, N)$
\[  \sup_{M\geq m_2} \left|\frac{1}{NM} \log \mu_N\left(P_{N,M}(x,y)\right) + \frac{\log Z_N}{N}\right| \leq \delta. \]
\end{proof}

Recall that for $(x, y) = (x_m, y_m)_{m\in \mathbb{N}}$ in $\left(\Omega|_N\right)^{\mathbb{N}}$
we have denoted 
\[  Q_{N,M}(x,y) = \left\{\left(\sum_{m=1}^\infty \frac{x^\prime_m}{a^m}, \sum_{m=1}^\infty \frac{y^\prime_m}{b^m}\right) 
\in X_\Omega |_N \middle|\, 
    \parbox{2.5in}{\centering $(x^\prime_m, y^\prime_m) \in \Omega|_N$ for all $m\geq 1$ with \\
    $x^\prime_m = x_m \> (1\leq m \leq \lfloor wM \rfloor)$ and \\ $y^\prime_m = y_m \> (1\leq m \leq M)$} \right\}.  \]
This is the image of the set $P_{N,M}(x,y)$ under the map $\left(\Omega|_N\right)^{\mathbb{N}}\to X_\Omega|_N$.
It is straightforward to check that for $(x,y) = (x_m, y_m)_{m\in \mathbb{N}}$ and 
$(x^\prime,y^\prime) = (x^\prime_m, y^\prime_m)_{m\in \mathbb{N}}$ in $\left(\Omega|_N\right)^{\mathbb{N}}$
the following three conditions are equivalent to each other:
\begin{itemize}
   \item $P_{N,M}(x,y) = P_{N,M}(x^\prime, y^\prime)$.
   \item $Q_{N,M}(x,y) = Q_{N,M}(x^\prime, y^\prime)$.
   \item $\left(x_1, \dots, x_{\lfloor wM\rfloor}, y_1, \dots, y_M\right) = 
   \left(x^\prime_1, \dots, x^\prime_{\lfloor wM\rfloor}, y^\prime_1,\dots, y^\prime_M\right)$.
\end{itemize}

For subsets $E$ and $F$ of $X_\Omega|_N$ we set 
\[  \dist_\infty(E,F) = \inf_{\substack{x\in E \\ y\in F}} \norm{x-y}_\infty. \]
Here $\norm{x-y}_\infty$ is the $\ell^\infty$-distance on $X_\Omega|_N\subset \mathbb{R}^{2N}$.

\begin{lemma} \label{lemma: estimate of duplication}
Let $N$ and $M$ be natural numbers.
Let $\left(x^{(1)},y^{(1)}\right), \dots, \left(x^{(k)}, y^{(k)}\right) \in \left(\Omega|_N\right)^{\mathbb{N}}$ with $k\geq 4^N+1$.
Suppose that 
\[   P_{N,M}\left(x^{(i)}, y^{(i)}\right) \neq P_{N,M}\left(x^{(j)},y^{(j)}\right) \quad \text{for $i\neq j$}. \]
Then there are $i$ and $j$ for which
\[  \mathrm{dist}_\infty \left(Q_{N,M}\left(x^{(i)},y^{(i)}\right), Q_{N,M}\left(x^{(j)}, y^{(j)}\right)\right) \geq b^{-M}. \]
\end{lemma}

\begin{proof}
The assumption $P_{N,M}\left(x^{(i)}, y^{(i)}\right) \neq P_{N,M}\left(x^{(j)},y^{(j)}\right)$ implies
\[  \left(x^{(i)}_1, \dots, x^{(i)}_{\lfloor wM\rfloor}, y^{(i)}_1, \dots, y^{(i)}_M\right) 
\neq \left(x^{(j)}_1, \dots, x^{(j)}_{\lfloor wM\rfloor}, y^{(j)}_1, \dots, y^{(j)}_M\right). \]
Since $k\geq 4^N+1$, there are $i$ and $j$ for which we have either 
\begin{equation}  \label{eq: x coordinates are distant}
      \norm{\sum_{m=1}^{\lfloor wM\rfloor} \frac{x_m^{(i)}}{a^m} - 
      \sum_{m=1}^{\lfloor wM\rfloor} \frac{x^{(j)}_m}{a^m}}_\infty \geq   2a^{-\lfloor wM\rfloor}
\end{equation}
or 
\begin{equation} \label{eq: y coordinates are distant}
   \norm{\sum_{m=1}^M \frac{y_m^{(i)}}{b^m} - \sum_{m=1}^M \frac{y^{(j)}_m}{b^m}}_\infty \geq   2b^{-M}.
\end{equation}
We have 
\[   Q_{N,M}(x,y) \subset \left(\sum_{m=1}^{\lfloor wM\rfloor} \frac{x_m}{a^m}, \sum_{m=1}^M \frac{y_m}{b^m}\right) + 
      \left[0, a^{-\lfloor wM\rfloor}\right]^N \times \left[0, b^{-M}\right]^N. \]
Therefore (\ref{eq: x coordinates are distant}) implies 
\[  \dist_\infty \left(Q_{N,M}\left(x^{(i)},y^{(i)}\right), Q_{N,M}\left(x^{(j)}, y^{(j)}\right)\right) 
\geq a^{-\lfloor wM\rfloor} \geq b^{-M}. \]
The condition (\ref{eq: y coordinates are distant}) implies 
\[   \dist_\infty \left(Q_{N,M}\left(x^{(i)},y^{(i)}\right), Q_{N,M}\left(x^{(j)}, y^{(j)}\right)\right) \geq  b^{-M}. \]
So we get the statement in both of the cases.
\end{proof}

\begin{proposition}   \label{prop: upper bound on epsilon Hausdorff dimension of projections of carpets}
For any $\delta>0$ there exists $\varepsilon>0$ such that for any $N\geq 1$
\[   \dimh\left(X_\Omega|_N, \norm{\cdot}_\infty, \varepsilon\right) \geq \log_b Z_N - \delta N. \]
\end{proposition}

\begin{proof}
If $\Omega|_N$ is a single point, then $Z_N = 1$ and the statement holds trivially.
So we assume that $\Omega|_N$ contains at least two points.

Since $Z_N$ is sub-multiplicative in $N$, there exists a positive constant $s$ satisfying 
$\log_b Z_N \leq s N$ for all $N\geq 1$.
Let 
\[   m_2 = m_2\left(\frac{\delta \log b}{2} \right)  \]
be the natural number introduced in Lemma \ref{lemma: upper bound on measure of P_NM}.
Let $N$ be a natural number.
By Lemma \ref{lemma: upper bound on measure of P_NM} there exists a Borel subset 
$R = R\left(\frac{\delta\log b}{2}, N\right)\subset \left(\Omega|_N\right)^{\mathbb{N}}$ satisfying $\mu_N(R) \geq 1/2$ and 
for any $M \geq m_2$ and $(x,y)\in \left(\Omega|_N\right)^{\mathbb{N}}$
\[  \left|\frac{1}{NM} \log \mu_N\left(P_{N,M}(x,y)\right) + \frac{\log Z_N}{N}\right| \leq \frac{\delta\log b}{2}. \]
Take $\varepsilon>0$ satisfying 
\begin{equation}  \label{eq: choice of epsilon for the lower bound on dimension of carpets}
   \varepsilon < b^{-m_2-1}, \quad 4 b^s \varepsilon^{\delta/2} < \frac{1}{2}. 
\end{equation}   

Suppose we are given a covering $X_\Omega|_N = \bigcup_{k=1}^\infty E_k$ with 
$0< \diam\left(E_k, \norm{\cdot}_\infty \right) < \varepsilon$ for all $k\geq 1$.
We set $D(E_k) := \diam\left(E_k, \norm{\cdot}_\infty \right)$ for simplicity of the notation.
We will show that 
\[   \sum_{k=1}^\infty D(E_k)^{\log Z_N - \delta N} > 1. \]
This proves $\dimh\left(X_\Omega|_N, \norm{\cdot}_\infty, \varepsilon\right) \geq \log Z_N - \delta N$.

For each $E_k$ we take a natural number $M_k \geq m_2$ satisfying 
\[  b^{-M_k-1} \leq D(E_k) < b^{-M_k}. \]
We define 
\[  \mathcal{C}_k = \left\{P_{N,M_k}(x,y)\middle|\, 
     (x,y)\in R \text{ with } Q_{N,M_k}(x,y)\cap E_k \neq \emptyset \right\}. \]
For $P_{N,M_k}(x,y)$ and $P_{N,M_k}(x^\prime,y^\prime)$ in $\mathcal{C}_k$, we have 
\[  \dist_\infty\left(Q_{N,M_k}(x,y), Q_{N,M_k}(x^\prime, y^\prime)\right) \leq  D(E_k) < \varepsilon < b^{-M_k}. \]
By Lemma \ref{lemma: estimate of duplication} we have $\left|\mathcal{C}_k\right| \leq 4^N$.

We have 
\[   R\subset \bigcup_{k=1}^\infty \bigcup_{P\in \mathcal{C}_k} P. \]
Hence 
\begin{equation} \label{eq: measure of R is one half}
     \frac{1}{2} \leq \mu_N(R) \leq \sum_{k=1}^\infty \sum_{P\in \mathcal{C}_k} \mu_N(P). 
\end{equation}     
Since $M_k \geq m_2$, every $P = P_{N,M_k}(x,y) \in \mathcal{C}_k$ with $(x,y)\in R$ satisfies 
\[   \frac{1}{N M_k} \log \mu_N(P)  \leq - \frac{\log Z_N}{N} + \frac{\delta\log b}{2}. \]
Hence for $P\in \mathcal{C}_k$
\begin{align*}
    \mu_N(P) & \leq \exp \left\{M_k \left(- \log Z_N + \frac{N \delta \log b}{2}\right)\right\} \\
    & = \exp\left\{M_k \log b \left(-\log_b Z_N + \frac{\delta N}{2}\right)\right\} \\
    & = b^{-M_k \left(\log_b Z_N - \frac{\delta N}{2}\right)}.
\end{align*}
From $b^{-M_k -1} \leq D(E_k)$ we have $b^{-M_k} \leq b \cdot D(E_k)$. So 
\[  \mu_N(P) \leq \left(b\cdot D(E_k)\right)^{\log_b Z_N - \frac{\delta N}{2}}, \quad 
     (P\in \mathcal{C}_k). \]
Therefore by (\ref{eq: measure of R is one half})
\begin{align*}
   \frac{1}{2} & \leq \sum_{k=1}^\infty \sum_{P\in \mathcal{C}_k} \left(b\cdot D(E_k)\right)^{\log_b Z_N - \frac{\delta N}{2}} \\
   & \leq  \sum_{k=1}^\infty 4^N \left(b\cdot D(E_k)\right)^{\log_b Z_N - \frac{\delta N}{2}}  \quad 
   (\text{by $\left|\mathcal{C}_k\right|  \leq 4^N$}) \\
   &< \sum_{k=1}^\infty 4^N b^{\log _b Z_N} D(E_k)^{ \log_b Z_N - \frac{\delta N}{2}} \\
   &\leq \sum_{k=1}^\infty 4^N b^{sN} D(E_k)^{\log_b Z_N - \frac{\delta N}{2}}  \quad (\text{by $\log_b Z_N \leq sN$})\\
   &= \sum_{k=1}^\infty 4^N b^{sN} D(E_k)^{\frac{\delta N}{2}} \cdot D(E_k)^{\log_b Z_N - \delta N} \\
   &< \sum_{k=1}^\infty 4^N b^{sN} \varepsilon^{\frac{\delta N}{2}} \cdot D(E_k)^{\log_b Z_N - \delta N} \quad 
       (\text{by $D(E_k) < \varepsilon$}) \\
   &< \sum_{k=1}^\infty \frac{1}{2^N} \cdot D(E_k)^{\log_b Z_N - \delta N} \quad 
      \left(\text{by $4 b^s \varepsilon^{\delta/2} < \frac{1}{2}$
      in (\ref{eq: choice of epsilon for the lower bound on dimension of carpets})}\right).
\end{align*}
Thus 
\[  1 <  \sum_{k=1}^\infty  D(E_k)^{\log_b Z_N - \delta N}. \]
\end{proof}

\begin{corollary} \label{cor: lower bound on mean Hausdorff dimension of carpets}
\[    \lmdimh\left(X_\Omega, \sigma, d\right) \geq \frac{\htop^w\left(\pi, \Omega, \sigma\right)}{\log b}. \]
\end{corollary}

\begin{proof}
For any natural number $N$,
the natural projection $(X_\Omega, d_N) \to (X_\Omega|_N, \norm{\cdot}_\infty)$ is one-Lipschitz.
So for any $\varepsilon>0$
\[  \dimh\left(X_\Omega, d_N, \varepsilon\right) \geq \dimh\left(X_\Omega|_N, \norm{\cdot}_\infty, \varepsilon \right). \]
By Proposition \ref{prop: upper bound on epsilon Hausdorff dimension of projections of carpets}, for any $\delta>0$ there exists
$\varepsilon>0$ so that for any $N\geq 1$
\[  \dimh\left(X_\Omega, d_N, \varepsilon\right) \geq \dimh\left(X_\Omega|_N, \norm{\cdot}_\infty, \varepsilon \right)
     \geq \log_b Z_N - \delta N. \]
Hence 
\begin{align*}
  \liminf_{N\to \infty} \frac{\dimh\left(X_\Omega, d_N, \varepsilon\right)}{N} & \geq 
     \lim_{N\to \infty} \frac{\log_b Z_N}{N} - \delta = \frac{1}{\log b} \lim_{N\to \infty} \frac{\log Z_N}{N} - \delta \\
  &= \frac{\htop^w\left(\pi, \Omega, \sigma\right)}{\log b}  - \delta, \quad 
    \left(\text{by $\htop^w\left(\pi, \Omega, \sigma\right) = \lim_{N\to \infty} \frac{\log Z_N}{N}$}\right).
\end{align*}     
Letting $\varepsilon\to 0$ and $\delta\to 0$, we get the statement.
\end{proof}

\begin{proposition}  \label{prop: upper bound on epsilon Hausdorff dimension of carpets}
For any $\varepsilon>0$ there exists a natural number $N_0 = N_0(\varepsilon)$ such that for any $N\geq N_0$
\[  \dimh\left(X_\Omega, d_{N-N_0}, \varepsilon\right) \leq (1+\varepsilon) \left(\log_b Z_N + \varepsilon N\right). \]
\end{proposition}

\begin{proof}
We will use Lemma \ref{lemma: upper bound on epsilon Hausdorff dimension}.
As in the proof of Proposition \ref{prop: upper bound on epsilon Hausdorff dimension of projections of carpets}
we take a positive number $s$ satisfying $\log_b Z_N \leq sN$ for all $N\geq 1$.

Let $N$ be a natural number.
We define a probability measure $\nu_N$ on $X_\Omega|_N$ as the push-forward of $\mu_N$ under the map
\[   \left(\Omega|_N\right)^{\mathbb{N}}\to X_\Omega|_N, \quad 
      \left(x_m, y_m\right)_{m \in \mathbb{N}} \mapsto
       \left(\sum_{m=1}^\infty \frac{x_m}{a^m}, \sum_{m=1}^\infty \frac{y_m}{b^m}\right). \]
We set 
\[  X(N) = X_\Omega|_N\times [0,1]^2\times [0,1]^2 \times \dots \subset \left([0,1]^2\right)^{\mathbb{N}}, \quad 
     \left(\text{recall } X_\Omega|_N \subset \left([0,1]^2\right)^N\right). \]
We have $X_\Omega \subset X(N)$. 
Recall that the metric $d$ on $\left([0,1]^2\right)^{\mathbb{N}}$ is defined by 
\[  d\left(\left(x_n, y_n\right)_{n\in \mathbb{N}}, \left(x^\prime_n, y^\prime_n\right)_{n\in \mathbb{N}}\right)
   = \sum_{n=1}^\infty  2^{-n} \max\left(|x_n-x^\prime_n|, |y_n-y^\prime_n|\right). \]

We denote the Lebesgue measure on the square $[0,1]^2$ by $Leb$.
We define
\[  \nu_N^\prime := \nu_N\otimes Leb \otimes Leb \otimes Leb \otimes \cdots. \]
This is a probability measure on $X(N)$.
For natural numbers $N, M$ and $(x,y)\in \left(\Omega|_N\right)^{\mathbb{N}}$ we set
\[  Q^\prime_{N,M}(x,y) := Q_{N,M}(x,y) \times [0,1]^2 \times [0,1]^2 \times \cdots \subset X(N). \]
We have 
\[   \nu^\prime_N\left(Q^\prime_{N,M}(x,y)\right) = \nu_N\left(Q_{N,M}(x,y)\right) \geq \mu_N\left(P_{N,M}(x,y)\right). \]

Given $\varepsilon>0$, we take $\delta>0$ with $6\delta^\varepsilon < 1$.
We take a natural number $K$ satisfying 
\begin{equation} \label{eq: choice of K for upper bound on dimension of carpets}
    a b^{-K} < \frac{\delta}{12}, \quad 
    \left(2a\right)^{s + \varepsilon}   \leq b^{\frac{\varepsilon K}{2}}. 
\end{equation}    
Let $L = L\left(\frac{\varepsilon \log b}{2}, K\right)$ be the natural number introduced in Lemma \ref{lemma: lower bound on measure of P_NM}.
Then by Lemma \ref{lemma: lower bound on measure of P_NM} for any $N\geq 1$ and any $(x,y)\in \left(\Omega|_N\right)^{\mathbb{N}}$
there exists $M \in [K, L]$ satisfying 
\[  \frac{1}{NM} \log \mu_N\left(P_{N,M}(x,y)\right) \geq - \frac{\log Z_N}{N} - \frac{\varepsilon\log b}{2}. \]
Then 
\begin{align*}
     \nu^\prime_N\left(Q^\prime_{N,M}(x,y)\right) & \geq \mu_N\left(P_{N,M}(x,y)\right) \\
     &\geq \exp\left\{NM\left(- \frac{\log Z_N}{N} - \frac{\varepsilon\log b}{2}\right)\right\} \\
     &= \exp\left\{-M\log b \left(\log_b Z_N + \frac{\varepsilon N}{2}\right)\right\} \\
     & = \left(b^{-M}\right)^{\log_b Z_N + \frac{\varepsilon N}{2}}.
\end{align*}

\begin{claim} \label{claim: lower bound on nu^prime_N measure}
For any $N\geq 1$ and any $(x,y)\in \left(\Omega_N\right)^{\mathbb{N}}$ there exists $M \in [K, L]$ satisfying
\[   \nu^\prime_N\left(Q^\prime_{N,M}(x,y)\right)  \geq \left(2a b^{-M}\right)^{\log_b Z_N + \varepsilon N}.  \]
\end{claim}

\begin{proof}
By the above argument,
there exists $M\in [K,L]$ satisfying 
\[    \nu^\prime_N\left(Q^\prime_{N,M}(x,y)\right) \geq \left(b^{-M}\right)^{\log_b Z_N + \frac{\varepsilon N}{2}}. \]
So it is enough to prove 
\[   \left(b^{-M}\right)^{\log_b Z_N + \frac{\varepsilon N}{2}}  \geq  \left(2a b^{-M}\right)^{\log_b Z_N + \varepsilon N}.  \]
This is equivalent to 
\begin{equation} \label{eq: b Z_N a}
   b^{\frac{\varepsilon M}{2}} \geq (2a)^{\frac{\log_b Z_N}{N} + \varepsilon}. 
\end{equation}   
From $\frac{\log Z_N}{N}\leq s$ and the choice of $K$ in (\ref{eq: choice of K for upper bound on dimension of carpets}),
\[  (2a)^{\frac{\log_b Z_N}{N} + \varepsilon} \leq (2a)^{s+\varepsilon} \leq b^{\frac{\varepsilon K}{2}}. \]
Since $M\geq K$, we have the above (\ref{eq: b Z_N a}).
\end{proof}

We take a natural number $N_0$ satisfying 
\[  \sum_{n>N_0} 2^{-n} < a b^{-L}. \]
For $u = (u_n)_{n\in \mathbb{N}}$ with $u_n\in [0,1]^2$, we denote $u|_N = (u_1, \dots, u_N)\in \left([0,1]^2\right)^N$.
Then for $u, v\in \left([0,1]^2\right)^{\mathbb{N}}$ and $N > N_0$
\begin{equation}  \label{eq: d and ell^infty norm of projections}
   d_{N-N_0}(u, v) < \norm{u|_N- v|_N}_\infty + a b^{-L}. 
\end{equation}   

\begin{claim}  \label{claim: diameter of Q^prime_NM}
For any $(x,y)\in \left(\Omega|_N\right)^{\mathbb{N}}$ and any natural numbers $N,M$ with $N> N_0$ and $M\leq L$, 
we have 
\[    0 < \diam\left(Q^\prime_{N,M}(x,y), d_{N-N_0}\right) < 2 a b^{-M}. \] 
Notice that $2ab^{-M} < \delta/6$ because we have assumed $a b^{-K} < \delta/12$
in (\ref{eq: choice of K for upper bound on dimension of carpets}).
\end{claim}

\begin{proof}
Obviously $Q^\prime_{N,M}(x,y) = Q_{N,M}(x,y)\times [0,1]^2\times [0,1]^2\times \cdots$ is not a single point.
So its diameter is positive.
By the above (\ref{eq: d and ell^infty norm of projections}) 
\[  \diam\left(Q^\prime_{N,M}(x,y), d_{N-N_0}\right) < \diam \left(Q_{N,M}(x,y), \norm{\cdot}_\infty\right) + a b^{-L}. \]
As we saw in (\ref{eq: diameter of Q_NM}) in \S \ref{subsection: calculation of metric mean dimension of carpets},
we have $\diam\left(Q_{N,M}(x,y), \norm{\cdot}_\infty\right) \leq ab^{-M}$.
So 
\[ \diam\left(Q^\prime_{N,M}(x,y), d_{N-N_0}\right) < a b^{-M} + a b^{-L} \leq 2 a b^{-M}, \quad  (\text{by $M\leq L$}).  \]
\end{proof}

By Claims \ref{claim: lower bound on nu^prime_N measure} and \ref{claim: diameter of Q^prime_NM}
for any natural number $N> N_0$ and any $(x,y)\in \left(\Omega|_N\right)^{\mathbb{N}}$ there exists 
a natural number $M\in [K,L]$ such that 
\begin{align*}
   & 0 < \diam \left(Q^\prime_{N,M}(x,y), d_{N-N_0}\right) < \frac{\delta}{6}, \\
   &\nu^\prime_N\left(Q^\prime_{N,M}(x,y)\right)  > \left(\diam \left(Q^\prime_{N,M}(x,y), d_{N-N_0}\right)\right)^{\log_b Z_N + \varepsilon N}. 
\end{align*}
Recall that we have assumed $6 \delta^\varepsilon < 1$.
Applying Lemma \ref{lemma: upper bound on epsilon Hausdorff dimension} to $\left(X(N), d_{N-N_0}\right)$, we get 
\[  \dimh\left(X(N), d_{N-N_0}, \delta\right) \leq (1+\varepsilon)\left(\log_b Z_N + \varepsilon N\right), \quad 
     (N> N_0).  \]
Since $0<\delta<\varepsilon$ and $X_\Omega\subset X(N)$, we have
\[  \dimh\left(X_\Omega, d_{N-N_0}, \varepsilon\right) \leq (1+\varepsilon) \left(\log_b Z_N + \varepsilon N\right), \quad 
     (N>N_0). \]
\end{proof}

\begin{corollary} \label{cor: upper bound on mean Hausdorff dimension of carpets}
\[  \umdimh\left(X_\Omega, \sigma, d\right) \leq \frac{\htop^w\left(\pi, \Omega, \sigma\right)}{\log b}. \]
\end{corollary}

\begin{proof}
Recall 
\[  \frac{\htop^w\left(\pi,\Omega,\sigma\right)}{\log b} = \lim_{N\to \infty} \frac{\log_b Z_N}{N}. \]
Let $\varepsilon$ be any positive number and let $N_0 = N_0(\varepsilon)$ be the natural number given by
Proposition \ref{prop: upper bound on epsilon Hausdorff dimension of carpets}.
Then for any $N>N_0$
\[  \frac{\dimh\left(X_\Omega, d_{N-N_0},\varepsilon\right)}{N} \leq (1+\varepsilon) \left(\frac{\log_b Z_N}{N} + \varepsilon\right). \]
Letting $N\to \infty$, we get 
\begin{align*}
   \limsup_{N\to \infty} \frac{\dimh\left(X_\Omega, d_N,\varepsilon\right)}{N}   
    &\leq (1+\varepsilon) \left(\lim_{N\to \infty} \frac{\log_b Z_N}{N} + \varepsilon\right) \\
    & = (1+\varepsilon) \left(\frac{\htop^w\left(\pi, \Omega, \sigma\right)}{\log b} + \varepsilon\right). 
\end{align*}    
Letting $\varepsilon \to 0$, we get the statement of the corollary.
\end{proof}

By Corollaries \ref{cor: lower bound on mean Hausdorff dimension of carpets} and
\ref{cor: upper bound on mean Hausdorff dimension of carpets}, we conclude
\[  \mdimh\left(X_\Omega, \sigma, d\right) = \frac{\htop^w\left(\pi, \Omega, \sigma\right)}{\log b}, \quad 
     (w = \log_a b). \]
This completes the proof of Theorem \ref{theorem: mean Hausdorff dimension and metric mean dimension of carpet systems}.

\appendix

\section{Mean Hausdorff dimension of $\left\{0,1,\frac{1}{2},\frac{1}{3},\dots, \right\}^\mathbb{N}$}
\label{appendix: mean Hausdorff dimension of simple example}

The purpose of this Appendix is to prove the result stated in Example \ref{example: 0,1,1/2,1/3,dots}.
Let 
\[  K = \left\{\frac{1}{n}\middle|\, n\geq 1\right\} \cup \{0\} = \left\{0,1,\frac{1}{2},\frac{1}{3},\dots\right\}. \]
Let $K^{\mathbb{N}}$ be the one-sided full shift on the alphabet $K$ with the shift map $\sigma: K^{\mathbb{N}}\to K^{\mathbb{N}}$.
We define a metric $d$ on it by 
\[ d\left((x_n)_{n\in \mathbb{N}}, (y_n)_{n\in \mathbb{N}}\right) = \sum_{n=1}^\infty 2^{-n} |x_n-y_n|. \]
In Example \ref{example: 0,1,1/2,1/3,dots} we claimed that 
\[  \mdimh\left(K^{\mathbb{N}}, \sigma, d\right) = 0, \quad 
     \mdimm\left(K^{\mathbb{N}},\sigma, d\right) = \frac{1}{2}. \]
We omit the proof of $\mdimm\left(K^{\mathbb{N}},\sigma, d\right) = \frac{1}{2}$ because it is rather straightforward.
In this Appendix we explain the detailed proof of $\mdimh\left(K^{\mathbb{N}}, \sigma, d\right) = 0$.
The proof is based on Lemma \ref{lemma: upper bound on epsilon Hausdorff dimension} in 
\S \ref{subsection: preparations for infinite dimensional carpets}.
We restate its special version (letting $c=1$ on Lemma \ref{lemma: upper bound on epsilon Hausdorff dimension}) 
here for the convenience of readers.

\begin{lemma}[$\subset$ Lemma \ref{lemma: upper bound on epsilon Hausdorff dimension}]
\label{revisited lemma: upper bound on epsilon Hausdorff dimension}
Let $\varepsilon$ and $s$ be positive numbers with $\varepsilon < 1/6$.
Let $(X, d)$ be a compact metric space with a Borel probability measure $\mu$.
Suppose that for any $x\in X$ there exists a Borel subset $A\subset X$ with $x\in A$ satisfying 
   \[   0< \diam\, A < \frac{\varepsilon}{6}, \quad \mu(A) \geq \left(\diam\, A\right)^s. \]
Then $\dimh(X, d, \varepsilon) \leq  2s$.
\end{lemma}

Set $X= K^{\mathbb{N}}$.
We define a probability measure $\nu$ on $K$ by 
\[  \nu\left(\{u\}\right) = au^2 \quad (u\neq 0), \quad 
     \nu\left(\{0\}\right) = \frac{1}{2}, \]
where $a$ is a positive number satisfying 
\[  a\left(1+\frac{1}{2^2} + \frac{1}{3^2} + \cdots\right) = \frac{1}{2}. \]
Indeed we have $a = \frac{3}{\pi^2}$, but we do not need its precise value.
We only need to use $a<1$.
We define $\mu = \nu^{\otimes \mathbb{N}}$. 
This is a Borel probability measure on $X$. 

Let $\varepsilon$ be a positive number with $\varepsilon<\frac{1}{6}$, and let $m$ be a natural number.
We take a positive number $\delta = \delta(\varepsilon, m)$ satisfying 
\begin{equation} \label{eq: choice of delta in Appendix}
   \delta < \min\left(\frac{\varepsilon}{12}, \frac{a^m}{8}, \left(\frac{1}{2}\right)^{\frac{m}{3}+1}\right). 
\end{equation}
We take a natural number $L$ satisfying 
\begin{equation*}  
     \sum_{n >L} 2^{-n} < \delta^{m^m}. 
\end{equation*}

\begin{claim} \label{claim: estimate of measure in Appendix}
For any $x\in X$ and any natural number $N$ there exists a Borel subset $A\subset X$ with $x\in A$ satisfying 
\[  0< \diam\, (A, d_N) < \frac{\varepsilon}{6}, \quad \mu(A) \geq \left(\diam(A, d_N)\right)^{\frac{6}{m}(N+L)}. \]
\end{claim}

We assume this claim for the moment.
Then by Lemma \ref{revisited lemma: upper bound on epsilon Hausdorff dimension}, for any natural number $N$
\[  \dimh(X, d_N, \varepsilon) \leq \frac{12}{m}(N+L). \]
Since $L$ is independent of $N$, we have 
\[  \limsup_{N\to \infty} \frac{\dimh(X, d_N, \varepsilon)}{N} \leq \frac{12}{m}.  \]
Letting $\varepsilon \to 0$ and $m\to \infty$, we conclude 
\[  \umdimh(X, \sigma, d) = 0. \]
So the rest of the problem is to prove Claim \ref{claim: estimate of measure in Appendix}.

\begin{proof}[Proof of Claim \ref{claim: estimate of measure in Appendix}]
Take $x = (x_n)_{n\in \mathbb{N}}$ with $x_n\in K$.
For this $x$ we introduce a partition 
\[  [1,N+L]\cap \mathbb{N} = I_0\cup I_1 \cup I_2 \cup\dots \cup I_m \cup I_{m+1} \quad (\text{disjoint union}) \]
by 
\begin{align*}
  &  I_0 = \left\{n \middle|\, x_n >\delta\right\}, \quad 
  I_k = \left\{n \middle|\, \delta^{m^k} < x_n \leq \delta^{m^{k-1}}\right\} \quad (1\leq k \leq m), \\
  &  I_{m+1} = \left\{n \middle|\, x_n \leq \delta^{m^m}\right\}.
\end{align*} 
There exists $m_0\in \{0,1,2,\dots, m\}$ satisfying 
\begin{equation}  \label{eq: choice of m_0 in Appendix}
    |I_{m_0}| \leq \frac{N+L}{m+1} < \frac{N+L}{m}. 
\end{equation}    
Set 
\[  r = \delta^{m^{m_0}}. \]
We have $r \leq  \delta < \varepsilon/12$.
We define $A\subset X$ by 
\[  A = \left\{(y_n)_{n\in \mathbb{N}}\in K^{\mathbb{N}}\middle|\, x_n-r  \leq  y_n  \leq x_n  \text{ for all $1\leq n\leq N+L$}\right\}. \]
This is not a single point. So its diameter is positive.
We have $x\in A$ and 
\[  \diam(A, d_N) \leq r + \sum_{n>L} 2^{-n} < r + \delta^{m^m} \leq 2r < \frac{\varepsilon}{6}. \]
We need to estimate 
\[   \mu(A)  = \prod_{n=1}^{N+L} \nu\left([x_n-r, x_n]\right). \]
Let $1\leq n \leq N+L$.
  \begin{enumerate}
      \item[Case 1.]   If $n\in I_k$ with $0\leq k < m_0$ then $x_n > \delta^{m^k}$ and hence 
        \begin{align*}
            \nu\left([x_n-r,x_n]\right) &\geq \nu\left(\{x_n\}\right) = a x_n^2  \\
             & > a \delta^{2m^k} \\
             & \geq  a \delta^{2m^{m_0-1}}    
              = a \left(\delta^{m^{m_0}}\right)^{\frac{2}{m}} = a r^{\frac{2}{m}} \\
             & > (2r)^{\frac{3}{m}} \quad \text{by $r \leq \delta < \frac{a^m}{8}$ in (\ref{eq: choice of delta in Appendix})} \\
             & \geq \left(\diam(A, d_N)\right)^{\frac{3}{m}}.
        \end{align*}
      \item[Case 2.]   If $n\in I_{m_0}$ then $x_n > \delta^{m^{m_0}} = r$ and hence 
         \begin{align*}
            \nu\left([x_n-r,x_n]\right) &\geq \nu\left(\{x_n\}\right) = a x_n^2  \\
             & > a r^2  \\
             & > (2r)^3 \quad \text{by $r\leq \delta < \frac{a^m}{8} < \frac{a}{8}$} \\
             & \geq \left(\diam(A, d_N)\right)^3.
         \end{align*}
         Since we know that $|I_{m_0}| \leq (N+L)/m$ by (\ref{eq: choice of m_0 in Appendix}), we have
         \[  \prod_{n\in I_{m_0}} \nu\left([x_n-r,x_n]\right) \geq     \left(\diam(A, d_N)\right)^{\frac{3}{m}(N+L)}.  \]
      \item[Case 3.]  If $n \in I_k$ with $k > m_0$ then $x_n \leq \delta^{m^{k-1}} \leq \delta^{m^{m_0}} = r$. 
               This implies $0\in [x_n-r,x_n]$. Hence
         \begin{align*}
           \nu\left([x_n-r,x_n]\right) & \geq \nu\left(\{0\}\right) = \frac{1}{2} \\
           & > (2r)^{\frac{3}{m}}  \quad \text{by $r \leq \delta < \left(\frac{1}{2}\right)^{\frac{m}{3}+1}$} \\
           & \geq \left(\diam(A,d_N)\right)^{\frac{3}{m}}.         
         \end{align*}         
   \end{enumerate}
Summarizing the above, we get 
\[  \mu(A) \geq  \underbrace{\left(\diam(A, d_N)\right)^{\frac{3}{m}(N+L)}}_{\text{contributions from Cases 1 and 3}}
     \cdot \underbrace{\left(\diam(A, d_N)\right)^{\frac{3}{m}(N+L)}}_{\text{contribution from Case 2}}
     = \left(\diam(A, d_N)\right)^{\frac{6}{m}(N+L)} \]
\end{proof}


\begin{thebibliography}{99}



\bibitem[BF09]{Barral--Feng_arXiv}
J.~Barral,  D.-J. Feng,
Weighted thermodynamic formalism and applications,
arXiv:0909.4247.





\bibitem[BF12]{Barral--Feng}
J.~Barral, D.-J. Feng, 
Weighted thermodynamic formalism on subshifts and applications, 
Asian J. Math. \textbf{16} (2012) 319-352.



\bibitem[BP17]{Bishop--Peres}
C.~J.~Bishop, Y.~Peres,
Fractals in probability and analysis,
Cambridge Studies in Advanced Mathematics, 162. Cambridge University Press,
Cambridge 2017.




\bibitem[Bow73]{Bowen}
R.~Bowen, 
Topological entropy for noncompact subsets,
Trans. Amer. Math. Soc. \textbf{184} (1973) 125-136.



\bibitem[Bed84]{Bedford}
T.~Bedford, 
Crinkly curves, Markov partitions and dimension, 
Ph.D. Thesis, University of Warwick, 1984.



\bibitem[DZG98]{Dai--Zhou--Geng}
X.~Dai, Z.~Zhou, X.~Geng, Some relations between Hausdorff-dimensions 
and entropies, 
Sci. China Ser. A \textbf{41} (1998) 1068-1075.


\bibitem[Dur10]{Durrett}
R.~Durrett,
Probability: theory and examples,
Fourth edition, 
Cambridge Series in Statistical and Probabilistic Mathematics, 31.
Cambridge University Press, Cambridge, 2010.





\bibitem[EW11]{Einsiedler--Ward}
M.~Einsiedler, T.~Ward,
Ergodic theory with a view towards number theory,
Graduate Texts in Mathematics, 259. 
Springer-Verlag London, London, 2011.




\bibitem[Fal89]{Falconer}
K.~J.~Falconer,
Dimensions and measures of quasi self-similar sets,
Proc. Amer. Math. Soc. \textbf{106} (1989) 543-554.



\bibitem[Fel68]{Feller}
W.~Feller, 
An introduction to probability theory and its applications, Volume I,
Third edition, John Wiley \& Sons, Inc., New York-London-Sydney, 1968.





\bibitem[FH16]{Feng--Huang}
D.-J. Feng, W. Huang, 
Variational principle for weighted topological pressure,
J. Math. Pures Appl. \textbf{106} (2016) 411-452.






\bibitem[Fur67]{Furstenberg}
H.~Furstenberg,
Disjointness in ergodic theory, minimal sets, and a problem in Diophantine approximation,
Math. Systems Theory \textbf{1} (1967) 1-49.




\bibitem[Gro99]{Gromov}
M.~Gromov, 
Topological invariants of dynamical systems and spaces of holomorphic maps: I,
Math. Phys. Anal. Geom. \textbf{2} (1999) 323-415.



\bibitem[KP96]{Kenyon--Peres}
R.~Kenyon, Y.~Peres,
Hausdorff dimensions of sofic affine-invariant sets,
Israel J. Math. \textbf{94} (1996) 157-178.





\bibitem[Lin99]{Lindenstrauss}
E.~Lindenstrauss,
Mean dimension, small entropy factors and an embedding theorem,
Inst. Hautes \'{E}tudes Sci. Publ. Math. \textbf{89} (1999) 227-262.


 
 
\bibitem[LT19]{Lindenstrauss--Tsukamoto_double}
E.~Lindenstrauss, M.~Tsukamoto,
Double variational principle for mean dimension, 
Geom. Funct. Anal., \textbf{29} (2019) 1048-1109.
 
 


\bibitem[LW00]{Lindenstrauss--Weiss}
E.~Lindenstrauss, B.~Weiss,
Mean topological dimension,
Israel J. Math. \textbf{115} (2000) 1-24.



 
\bibitem[Mc84]{McMullen}
C.~McMullen,
The Hausdorff dimension of general Sierpinski carpets, 
Nagoya Math. J. \textbf{96} (1984) 1-9.
 
 

\bibitem[Mis04]{Misiurewicz}
M.~Misiurewicz,
On Bowen’s definition of topological entropy,
Discrete and Continuous Dynamical Systems \textbf{10}
(2004).




\bibitem[ST21]{Shinoda--Tsukamoto}
M.~Shinoda, M.~Tsukamoto,
Symbolic dynamics in mean dimension theory,
Ergod. Th. $\&$ Dynam. Sys.
\textbf{41} (2021) 2542-2560.



\bibitem[Tsu21]{Tsukamoto_weighted}
M.~Tsukamoto,
New approach to weighted topological entropy and pressure,
Ergod. Th. $\&$ Dynam. Sys. First view,
DOI: https://doi.org/10.1017/etds.2021.173.
arXiv:2108.03795.





\end{thebibliography}
\end{document}